\documentclass[12pt,a4,twoside]{amsart}
\pdfoutput=1
\usepackage{latexsym}
\usepackage{amssymb}
\usepackage{amsmath}
\usepackage{color}
\usepackage{graphicx}
\usepackage[all]{xy}
\newtheorem{theorem}{{\sc Theorem}}
\newtheorem{definition}{{\sc Definition}}
\newtheorem{lemma}{{\sc Lemma}}

\newtheorem{proposition}{{\sc Proposition}}

\newtheorem{remark}{{\sc Remark}}
\newtheorem{example}{{\sc Example}}

\newcommand{\RR}{{\mathbb R}}
\newcommand{\CC}{{\mathbb C}}

\def\colP#1#2{\{P_{#1},\cdots, P_{#2}\}}
\def\colQ#1#2{\{Q_{#1},\cdots, Q_{#2}\}}

\def\tP#1#2#3{\Delta(P_{#1}, P_{#2}, P_{#3})}
\def\P#1#2{\{P_{#1}, P_{#2}\}}
\def\PP#1#2#3{\{P_{#1}, P_{#2}, P_{#3}\}}
\def\PPP#1#2#3#4{\{P_{#1}, P_{#2}, P_{#3}, P_{#4}\}}

\newcommand{\CH}{\mathop{\mathrm{CH}}}
\newcommand{\Up}{\mathop{\mathrm{Up}}}
\newcommand{\Pow}{\mathop{\mathrm{Pow}}}
\newcommand{\Vor}{\mathop{\mathrm{Vor}}}
\newcommand{\Del}{\mathop{\mathrm{Del}}}
\newcommand{\Dom}{\mathop{\mathrm{Dom}}}
\newcommand{\Sep}{\mathop{\mathrm{Sep}}}
\newcommand{\Aff}{\mathop{\mathrm{Aff}}}
\newcommand{\Cone}{\mathop{\mathrm{Cone}}}
\newcommand{\Interior}{\mathop{\mathrm{Interior}}}
\newcommand{\sm}{\mathop{\mathrm{sm}}}
\newcommand{\grad}{\mathop{\mathrm{grad}}}
\newcommand{\Log}{\mathop{\mathrm{Log}}}

\begin{document}
\title{Power diagrams and Morse Theory}
\author{D.\ Siersma
  , M.\ van Manen}
 
\address{Mathematisch Instituut, Universiteit Utrecht, PO Box 80010, \ 3508 TA Utrecht The Netherlands.
\texttt {D.Siersma@uu.nl}  \and
Department of Mathematics, Hokkaido University
Kita 10, Nishi 8, Kita-Ku, Sapporo, Hokkaido, 060-0810, Japan.
\texttt{manen@math.sci.hokudai.ac.jp}}
\begin{abstract}
We study Morse theory of the (power) distance function to a set of points in $\RR^n$.
We describe the topology of the union of the corresponding set of growing balls by a Morse poset.
The Morse poset is related to the power tesselation of $\RR^n$.
We remark that the power diagrams from computer science are the spines of amoebas
in algebraic geometry, or the hypersurfaces in tropical geometry.
We show that there exists a discrete Morse function on the coherent triangulation,
dual to the power diagram, such that its critical set
equals the Morse poset of  the power diagram.
\end{abstract}
\maketitle
\setcounter{section}{0}
\section{Introduction}
\begin{figure}[h]
\begin{center}
  \includegraphics[height=4cm]{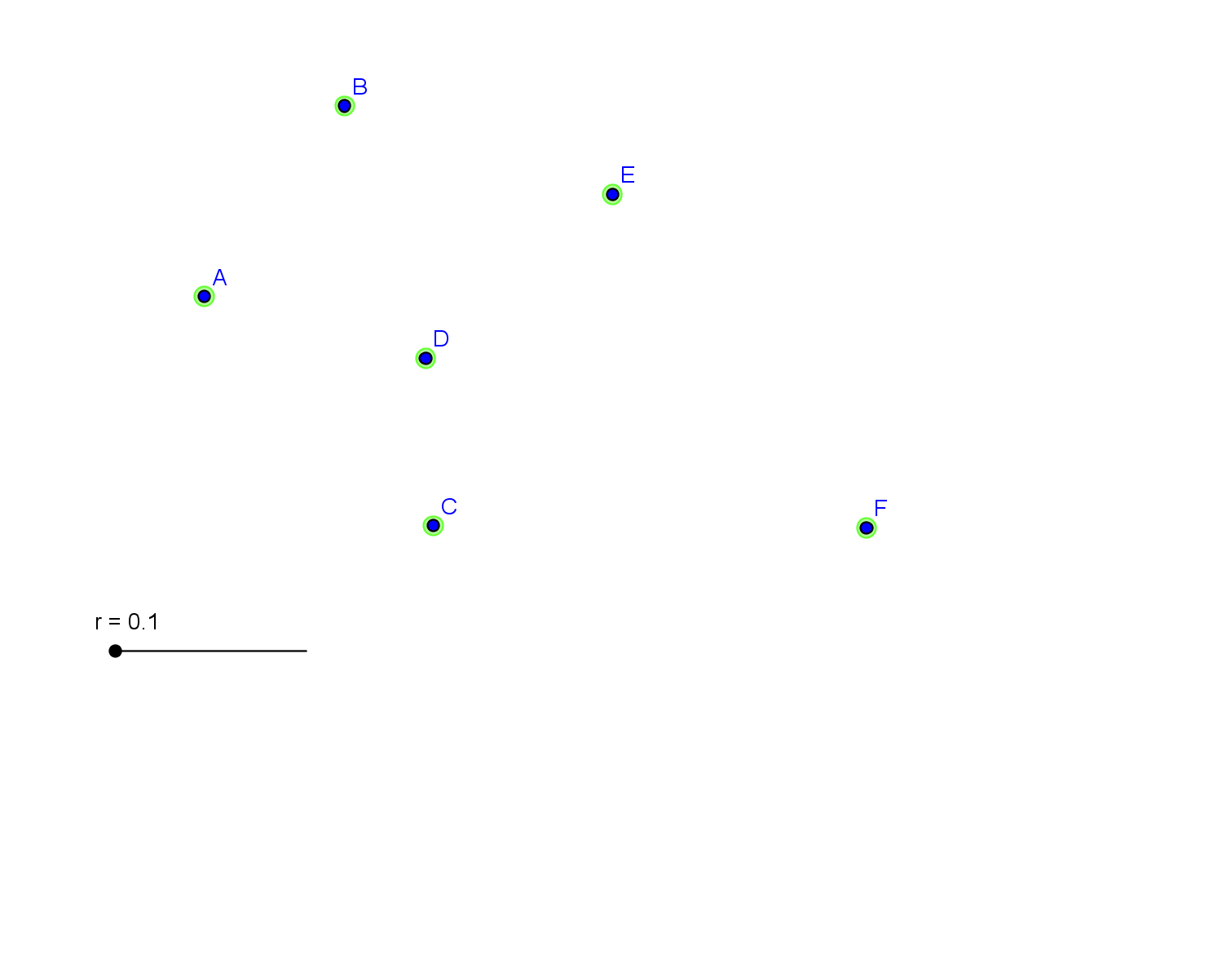}
\includegraphics[height=4cm]{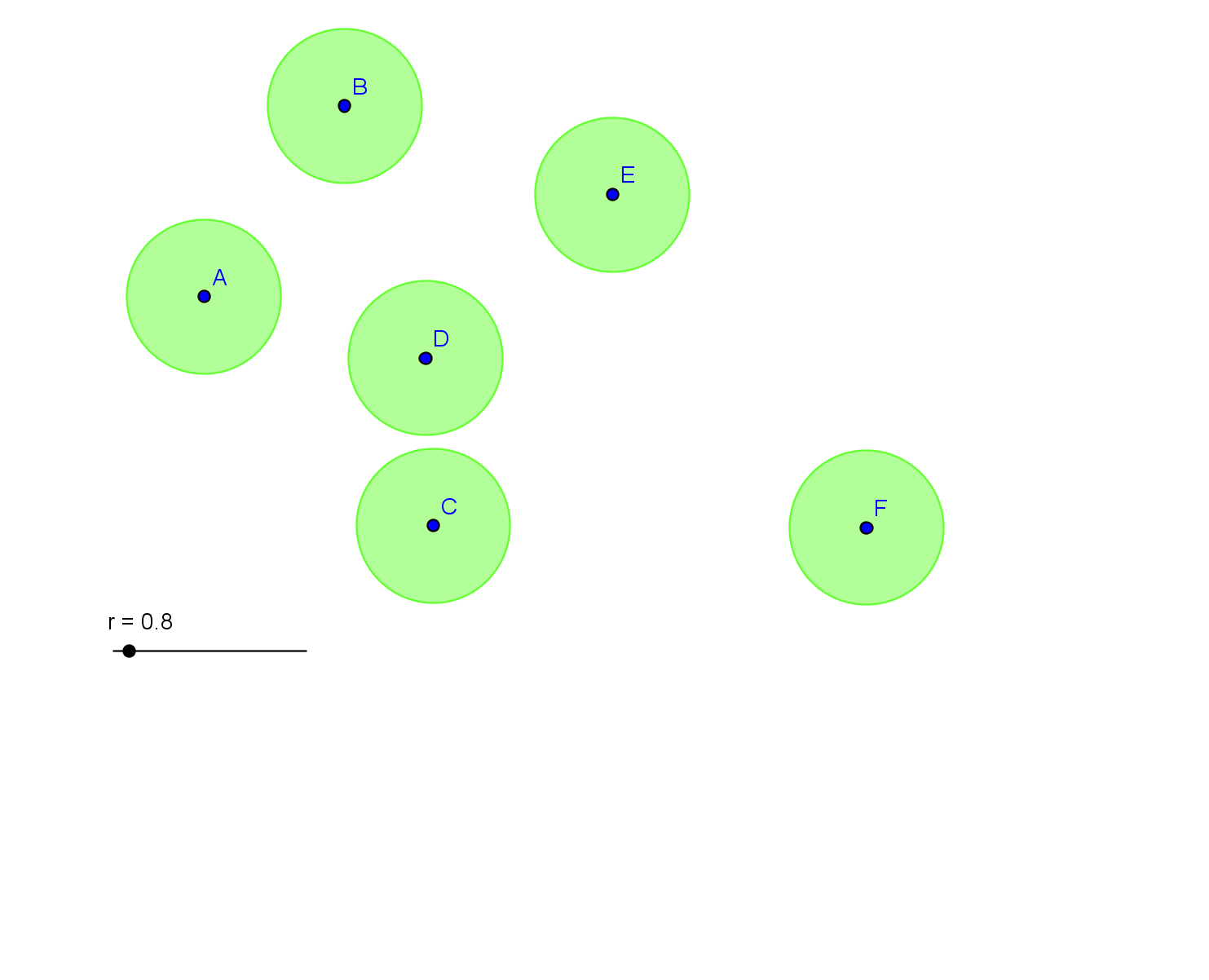}
\includegraphics[height=4cm]{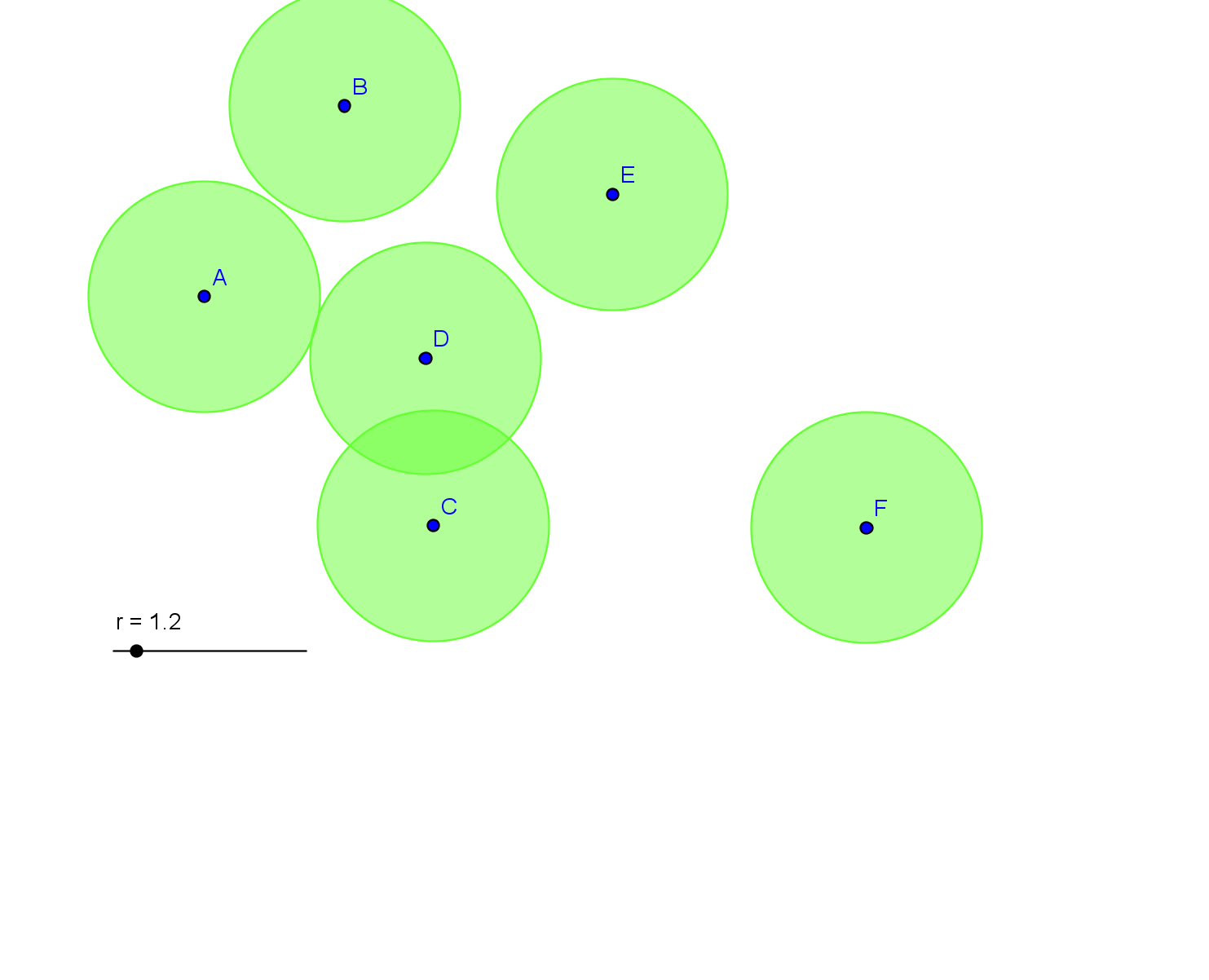}
\includegraphics[height=4cm]{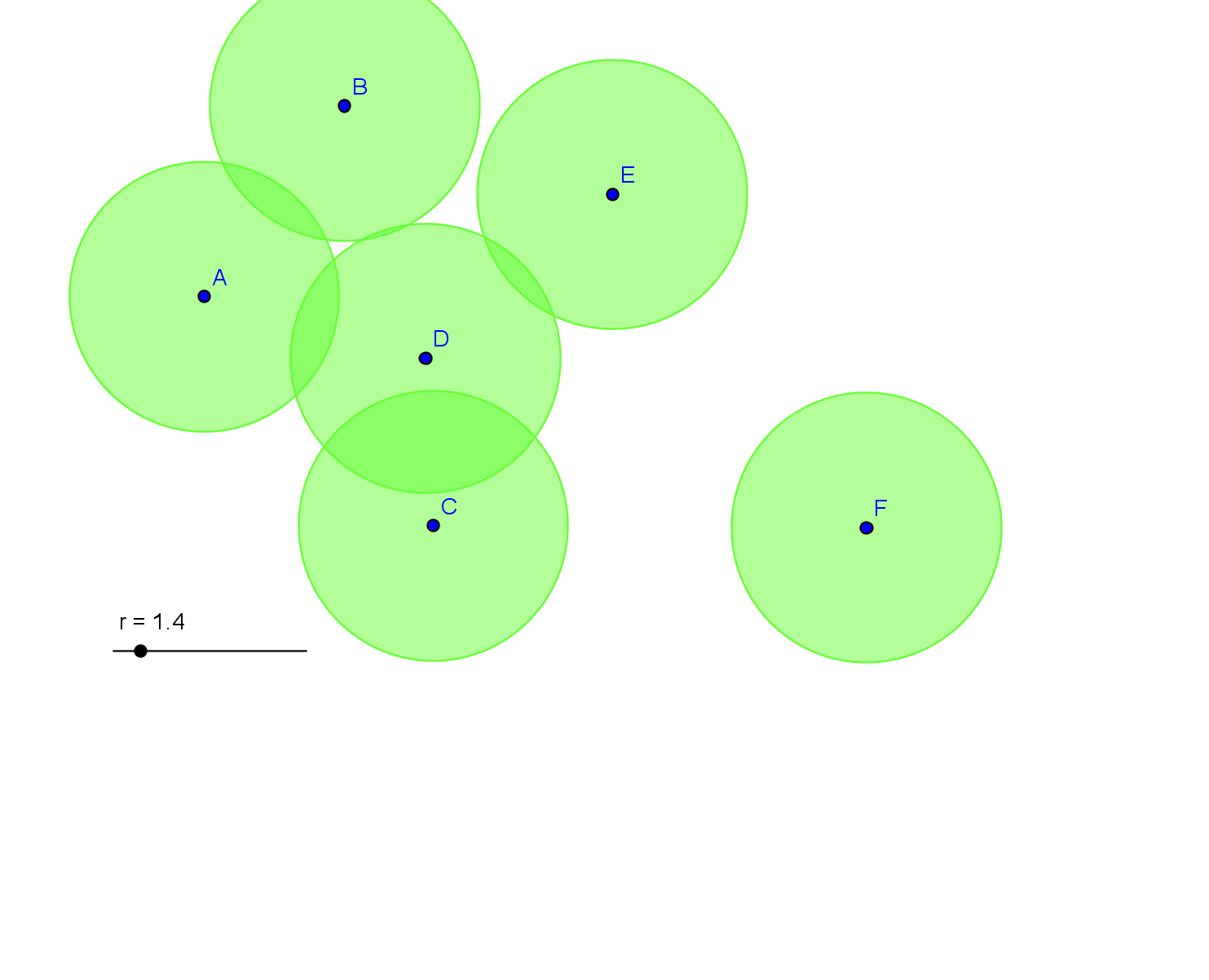}
\includegraphics[height=4cm]{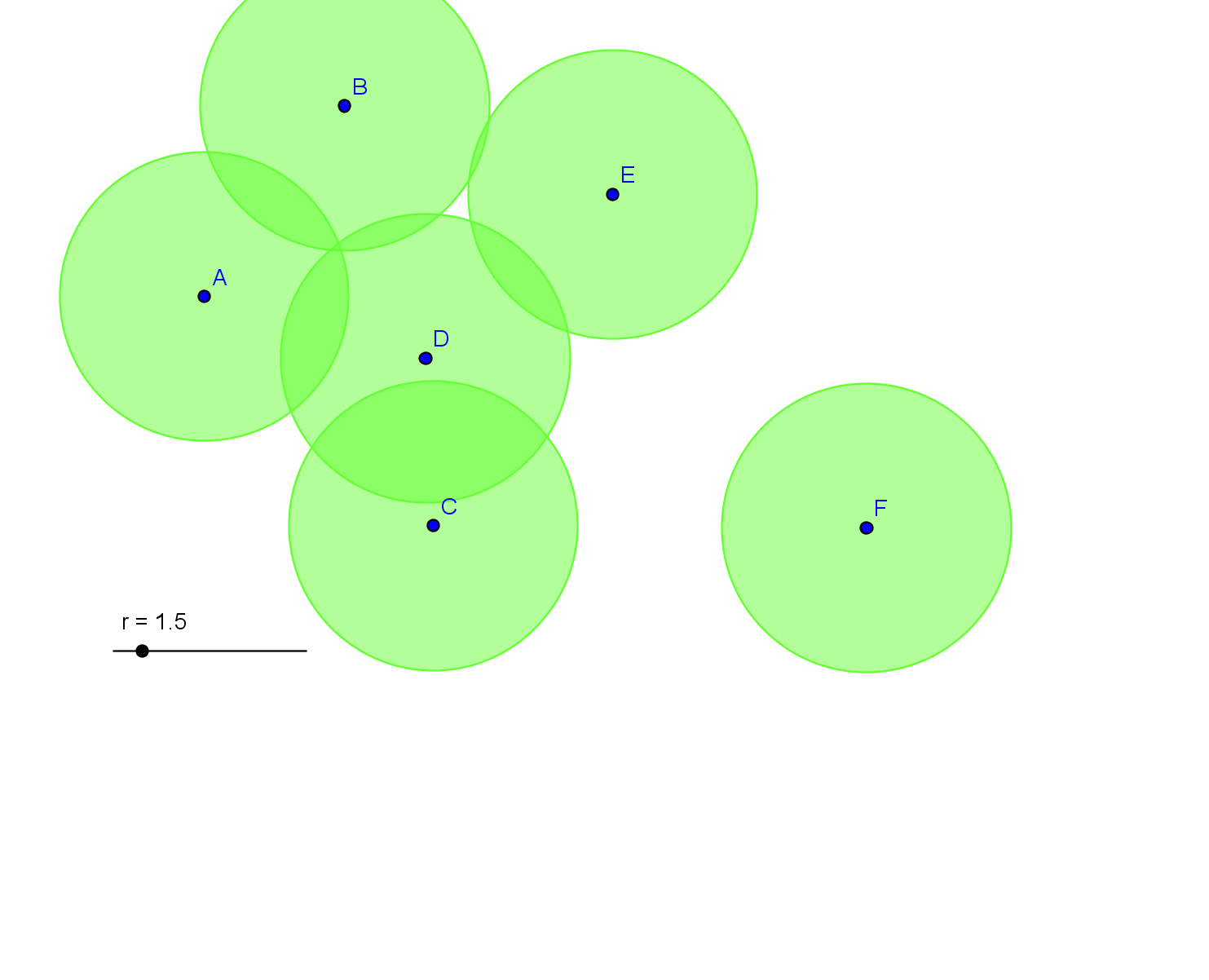}
\includegraphics[height=4cm]{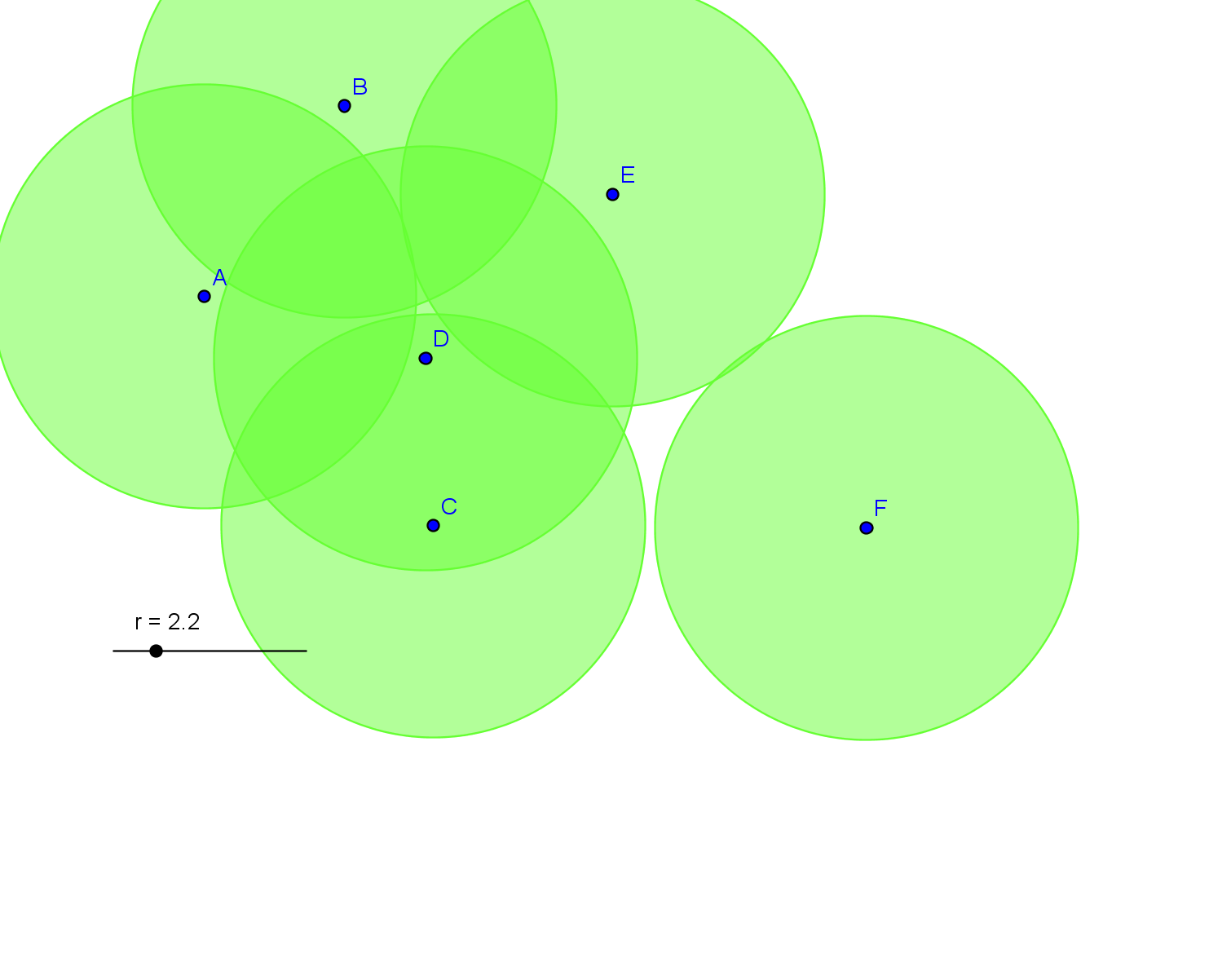}
\includegraphics[height=4cm]{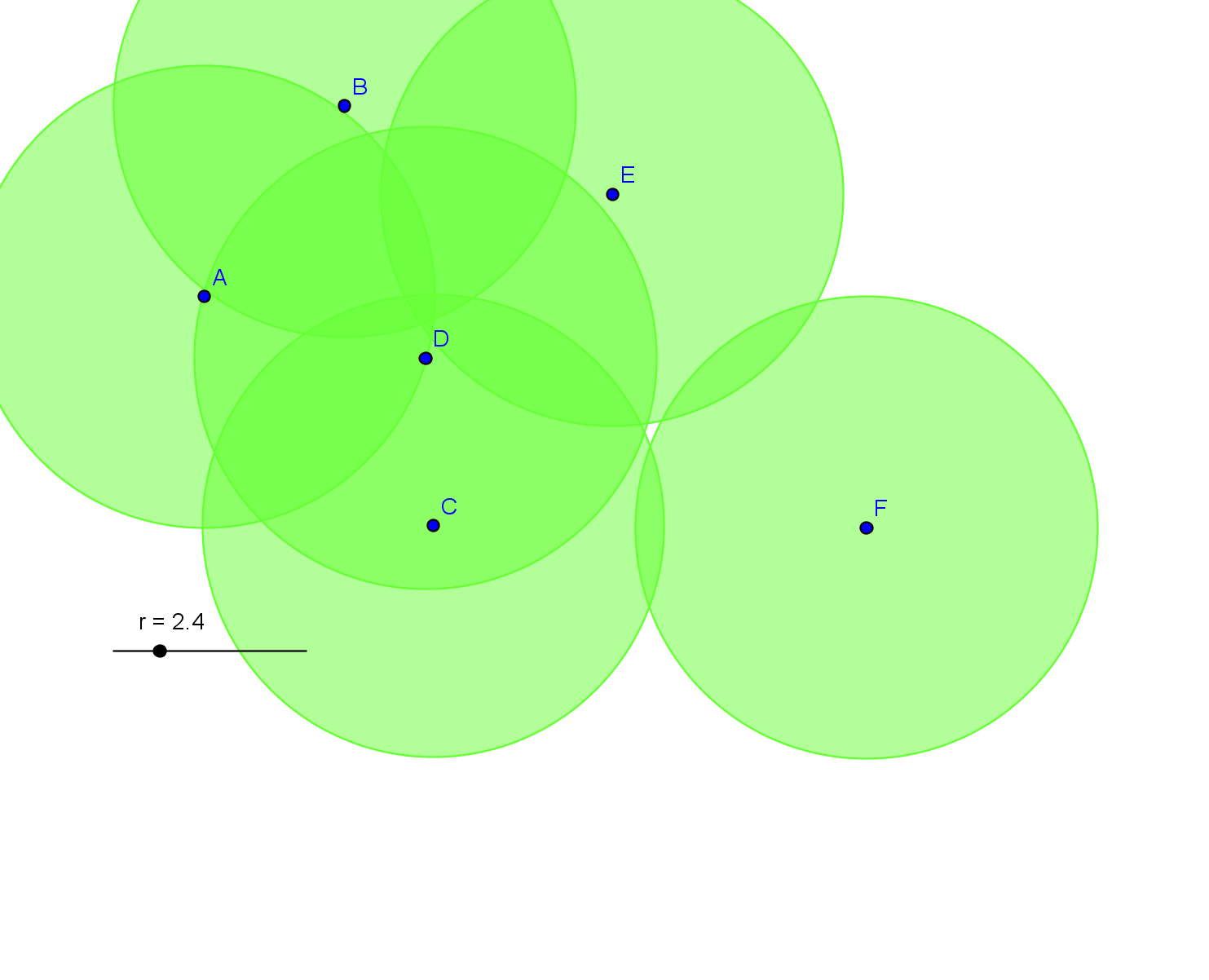}
\includegraphics[height=4cm]{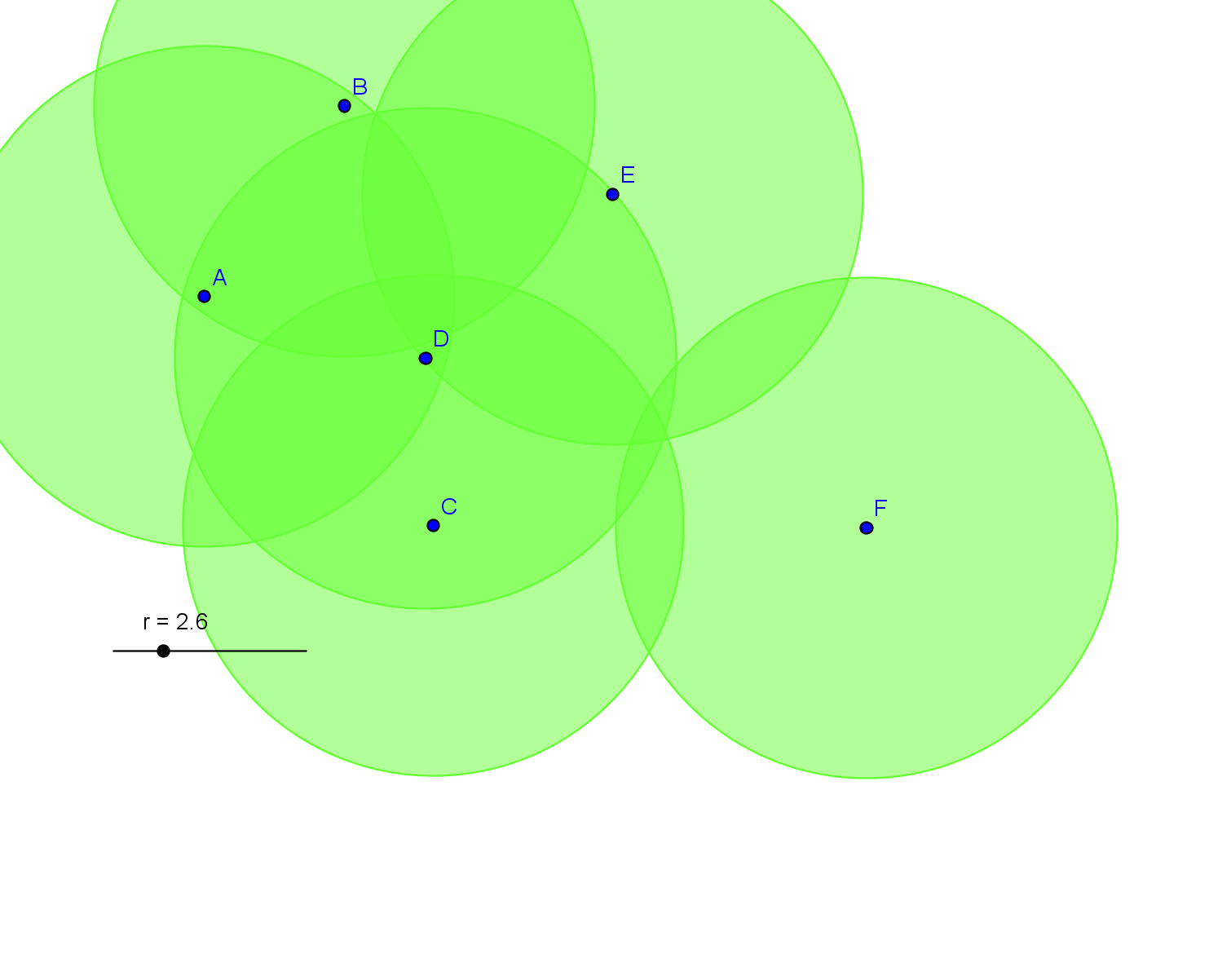}

\end{center}
\caption{Growing balls.} \label{fig:vm6}
\end{figure}

In this paper we study in a Euclidean space of arbitrary dimension the evolution of a set of balls with fixed centers $P_1,\cdots,P_N$ with respect to increasing radius. Let $B(P_i,r_i)$ be the ball with center $P_i$ and radius $r_i$. We let all balls grow simulaniously with increasing radius $r_i=t$. One can also see this as a set of homogeneous wave fronts that start from the centers with equal (and homogeneous) speed. See figure \ref{fig:vm6}.
 The space covered by the balls $B(P_1,t),\cdots , B(P_N,t)$ is the set $B(t) = 
\bigcup B(P_i,t)$. 
We will discuss several aspects of this process:
\begin{itemize}
\item[[T]] (T-Properties) change of topology B(t) if t increases.
\item[[G]] (G-Properties) the geometry of the set of points, which are covered for the first time, but exactly by two different balls; to
be more precise we get
\begin{equation*}
\Vor = \{ x \in \RR^n \,\mid\,\exists i \ne j\,: \,\,{}d(x,P_i) = d(x,P_j) \le d(x,P_k) \ \mathrm{for all} \ k \}.
\end{equation*}
This is the well-known {\em Voronoi diagram} of
the Euclidean distance function. Here $ d(x,P_i) = || x- P_i||$ is the Euclidean distance from $x$ to $P_i$ 
\end{itemize}

For the study of T-Properties we will consider Euler characteristic, topological Morse theory,
simplicial and cell-complexes and discrete Morse Theory.
G-Properties will lead us to Voronoi diagrams, Delaunay tesselations, and the
generalizations: power diagrams and power Delaunay tesselations.
\par

There is a second way to describe the union of balls $B(t)$. It can be done by the lower level sets of the minimum of the distance functions from a point $x$ in the Euclidean space to the centers $P_i$:
$$ d(x) = \min_{i=1,\cdots,N} || x- P_i||. $$
With this notation 
$B(t) = \{ x \in \RR^n | d(x) \le t \}  $.

In order to have differentiability in the points of $\colP1N$
and to have a nicer formula for the gradient, we study
\[ g(x) =   \min_{i=1,\cdots,N}  \frac12 || x- P_i||^2 = \min_{i=1,\cdots,N} g_i(x) \  \mbox{\rm where} \ \ g_i(x) =\frac12 || x- P_i||^2 ,\]
which behaves similarly to $d$, e.g. $d$ and $g$  have the same
set of level curves.
Note that $B(\sqrt{2t}) = \{ x \in \RR^n | g(x) \le t \} $.

The evolution of the set of balls will be studied with the (topological) Morse theory of the function $g: \RR^n \rightarrow \RR .$
The same scheme can be applied to other distance functions. From section \ref{sec-power} on we will use the
so-called power distance function, which generalizes the Euclidean distance.

The paper is organized as follows:
\par
We first describe in section \ref{sec-2d} a simple case : points in the plane with
Euclidean distance. We use elementary methods (as in \cite{VoronoiDist}) to
study the geometry and topology in order to give a good an elementary view on
what is going on. The treatment in later sections is more abstract and has greater
generality (power distance and any dimension).
\par
In section \ref{sec-power} we give the definition of power distance and review properties of the power diagram. 
This is a tesselation
of $\RR^n$ by means of the power distance, which an additively weighted 
distance function to a point set. A precise definition is given in section \ref{sec-power}.
\par
Power diagrams were introduced
as a generalization of Voronoi diagrams by Aurenhammer in \cite{MR873251}.
They play an important role in computational geometry.
We will show that power diagrams can also be defined in a second way: by affine functions and as such
they occur in algebraic geometry, as the spines of
amoebas, see \cite{MR2040284}. Another appearance of the power diagram is
what other authors call a tropical hypersurface, see \cite{AG0306366}.
They also appear when string theorists use toric geometry, see for instance
\cite{MR1969655}. We spell out the correspondence
between the classical theory of Voronoi and power diagrams and that of 
the spines of amoebas and tropical hypersurfaces. Also the dual of the power diagram, the Delaunay tesselation,
plays an important role in the rest of the paper. 
\par
In section \ref{morse} we study the critical points of the minimum of the power distances to the points $P_1,\ldots,P_N$ and  relate it to the Morse poset. This Morse poset is a partial ordered subset, or poset, of the Delaunay tesselation. It encodes the critical points and their Morse types.
It turns out that on the complement of the Morse posets there exist a discrete vector field in
the sense of Forman.
\par
We present in section \ref{sec-discr} the relation between power diagrams
and the discrete Morse theory introduced
by Forman in \cite{MR1612391}. This theory is important, because it describes the order of cell attaching
in a combinatorial model for the union of balls $B(t)$.
\par
In \cite{Wrapping} Edelsbrunner asked to ``elucidate the two approaches''. The one
being the study of embedded point sets using Euclidean distance functions, and the other one
being discrete Morse theory. We solve this in  section \ref{sec:morseposetToVf}, where we construct a discrete Morse function from Euclidean distance information.
\par
In applications one extensively uses the Morse theory point of view
to distinguish between different triangulations of point sets.
The Morse poset depends on the metric structure of the
point set and can be used as shape classifier. Examples were first given in \cite{SieVanM} and \cite{VoronoiDist}.
\par
Most of the results in this paper can in some form or another be found in the literature.
We tried to include appropriate references in the sections, but in all likelihood we are not complete.
New material is contained in section \ref{sec:morseposetToVf},  especially theorem \ref{sec:from-morse-poset} 
and \ref{sec:morse-ext}.  

\section{Two dimensional case, intuitive introduction}\label{sec-2d}
\subsection{Two dimensional case}
We start with a set of $N$ different points $\colP1N$
 in the plane $\RR^2$.
As mentioned above we study the function
\[ g(x) = \min_{i=1,\cdots,N}  \frac12 || x- P_i||^2 = \min_{i=1,\cdots,N} g_i(x) \  \mbox{\rm where} \ \ g_i(x) =\frac12 || x- P_i||^2 ,\]
Note that
\begin{equation*}
 \grad \ g_i (x) =  \grad \ \frac12 || x- P_i||^2  =  \overrightarrow{xP_i}
\end{equation*}
It follows that the set of points where $g$ is not differentiable
is exactly
\begin{equation*}
 \Vor = \{ x \in \RR^2 \,\mid\,\exists i \ne j\,: \,\,{}d(x,P_i) = d(x,P_j) \le d(x,P_k) \ \text{for all} \ k \}.
\end{equation*}
This is the well-known {\em Voronoi diagram} of
the Euclidean distance function $d$.
The {\em (closed) Voronoi cells} are defined by:
\begin{equation*}
\Vor (P_i)= \{ x \in \RR^2 \,\mid\,  d(x,P_i) \le d(x,P_k) \ \mbox{for all} \ k \}. 
\end{equation*}
The two dimensional Voronoi diagrams consist of
\begin{itemize}
\item 2-dimensional Voronoi cells,
\item 1-dimensional Vornonoi edges,
\item 0-dimensional Vornoi vertices.
\end{itemize}
The last two constitute a hypersurface, the Voronoi diagram.

For the theory of Voronoi diagrams we refer to Aurenhammer
\cite{Au}, Edelsbrunner \cite{Ed} and the book of
Okabe-Books-Sigihara \cite{OBS}. Voronoi diagrams have many
applications in mathematics and computer science, but also in
geography, biology, crystallography, marketing, cartography, etc.
\par
The level curves of the squared distance function $g$ can be considered as
{\em wave fronts}, which start from the points of $\colP1N$.
These wave fronts $\{ g = t \}$ bound regions $\{ g \le
t \}$, where the wave front has already passed, just as a
region passed by a forest fire. Each $\{ g \le
t \}$ is a union of balls.
\par
We want to study the change of topology of these regions $\{ g \le t \}$.
We start with an instructive and very simple example
with $3$ points, where two different positions of the points of
$\PP123$ give rise to different topological behavior. We first consider
the simplest indicator for topological changes: the Euler characteristic $\chi$.
\par
\begin{figure}[htb]
\begin{center}
 \includegraphics[width=\textwidth]{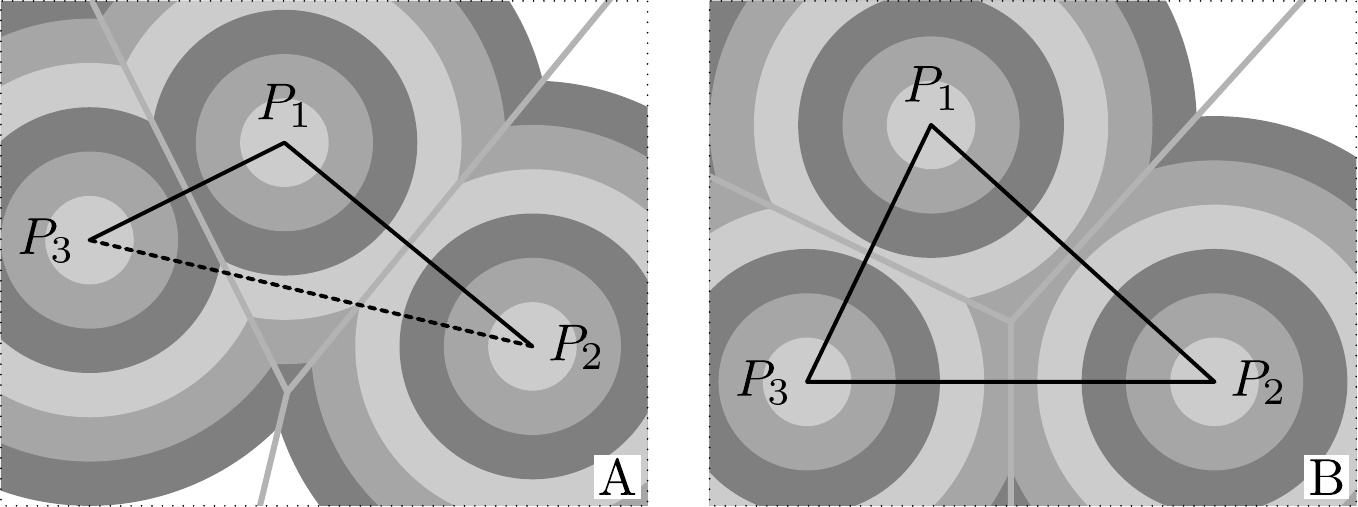}
\end{center}
\caption{Evolution of a wave front from three points, cases A and B}
\label{fig:fig0}
\end{figure}

Consider the level sets in figure \ref{fig:fig0}.
Initially the wave fronts surround three different regions.
So the Euler characteristic $\chi = 3$. We will report 
changes in $\chi$ as $t$ grows. Next, two regions meet in the center of $P_1P_2$ and we
get two contractible regions, so $\chi = 2$.  After that the third
region meets the other (combined) region in the center of $P_1P_3$, this gives $\chi = 1$.
\par
Next we distinguish two cases:
\par
In case A (figure \ref{fig:fig0}), where one of the angles is
{\it obtuse}, the region becomes bigger and bigger and $\chi$ does not
change anymore.
\par
In case B (figure \ref{fig:fig0})  all the angles are {\it acute} and
now the wave fronts meet once more: in the center of $P_2P_3$
 and enclose a region in the middle of the triangle.
The set $\{ g \le t \}$ forms a cycle and is no longer
contractible. Now $\chi = 0$.
If $t$ grows further then the enclosed region
in the middle disappears and this changes $\chi$ to $1$. A detailed picture is figure \ref{fig:fig4}B. 
There is
only one region left, which is contractible and there are no
changes if $t$ increases further.
\par

We intend to study this type of process in the present paper. The
special points where the wave fronts meet and the points where
they become non-differentiable are directly related to the Voronoi
diagram. As we see in the above example we need more refined
information than the Voronoi diagram in order to understand the
topological behavior of the wave fronts.
\subsection{Behavior of $g$ on the plane $\RR^2$.}\label{subs-crit}
The behavior of $g$ on the {\em interiors of the Voronoi cells}
is clear. In the points of $\colP1N$ the function $g$ has its
minimal value and there are no other special points in the
interiors of the Voronoi cells. The level curves are there smooth for $t > 0$.
Next we study the behavior of $g$ on neighborhoods of the points on the Voronoi
diagram: the edges and the vertices.
\par
 The edges of the Voronoi diagram are parts of perpendicular bisectors of two
points $P_i$ and $P_j$. For an example look at the points $P_2$ and $P_3$ in figure \ref{fig:fig0}.
\par
Let $\Vor(P_2,P_3) = \Vor (P_2) \cap \Vor (P_3)$ be the Voronoi edge between
the Voronoi cells of $P_2$ and $P_3$. We will call the perpendicular bisector
of the points $P_2$ and $P_3$ the separator
$\Sep(P_2,P_3)$ of $P_2$ and $P_3$.
Indeed, the perpendicular bisector separates the point $P_2$ from the point $P_3$.
Let $c_{23}=c(P_2,P_3)$ be the midpoint of the segment $P_2P_3$:   $c_{23} = \Sep(P_2,P_3) \cap P_2P_3$. 
\par
\begin{figure}[htb]
\begin{center}
 \includegraphics[width=\textwidth]{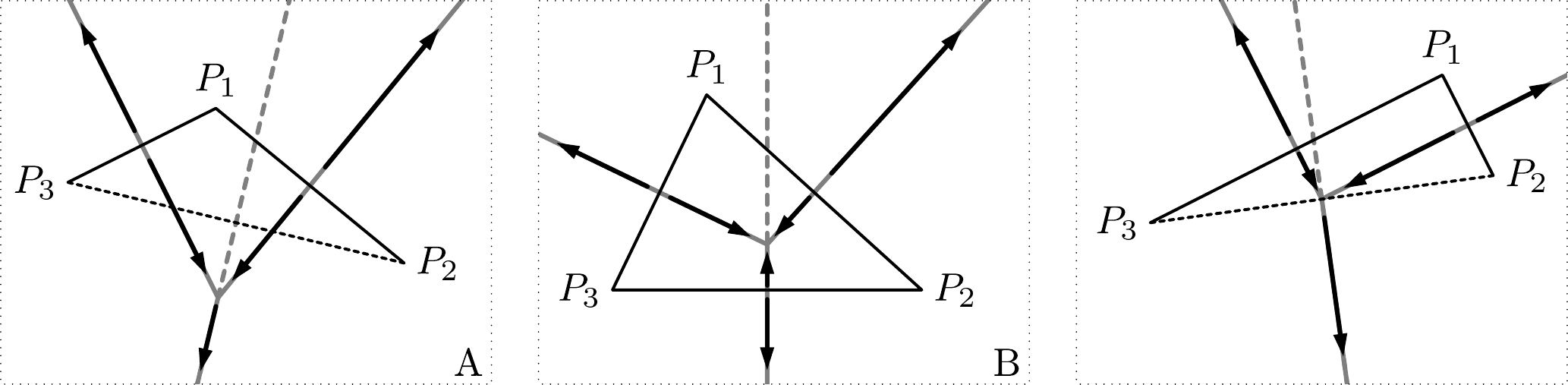}
\end{center}
\caption{Positions of $c_{23}$ with respect to $\Vor(P_2,P_3)$. The gray dashed line 
 is $\Sep(P_2,P_3)$. The arrows point to the direction in which $g$ increases on the
Voronoi diagram.}
\label{fig:fig1}
\end{figure}

The edge $\Vor(P_2,P_3)$ is contained in the separator: $\Vor(P_2,P_3)\subset \Sep(P_2,P_3)$.
The position of $c_{23}$ with respect to $\Vor(P_2,P_3)$ is important. There are three
cases (cf.\ figure \ref{fig:fig1}):
\begin{enumerate}
\item\label{item:2} $c_{23}$ lies outside $\Vor(P_2,P_3)$, the triangle $\tP123$ is obtuse. This is case A in figures \ref{fig:fig0} and \ref{fig:fig1}.
\item\label{item:3} $c_{23}$ lies in the interior of $\Vor(P_2,P_3)$, $\tP123$ is acute. This is case B in figures \ref{fig:fig0} and \ref{fig:fig1}.
\item\label{item:4} $c_{23}$ is a boundary point of $\Vor(P_2,P_3)$, $\tP123$ is right-angled. This is the third case in figure \ref{fig:fig1}.
\end{enumerate}
The position of $c_{23}$ with respect to $\Vor(P_2,P_3)$ determines the behavior of
$g$ on $\Vor(P_2,P_3)$.
In cases (~\ref{item:2}) and (~\ref{item:4}) $g$ is monotone on the edge; in case (~\ref{item:3}) $g$
is not monotone, but increasing from $c_{23}$ in both directions.
\par
\begin{figure}[htb]
\begin{center}
 \includegraphics[width=\textwidth]{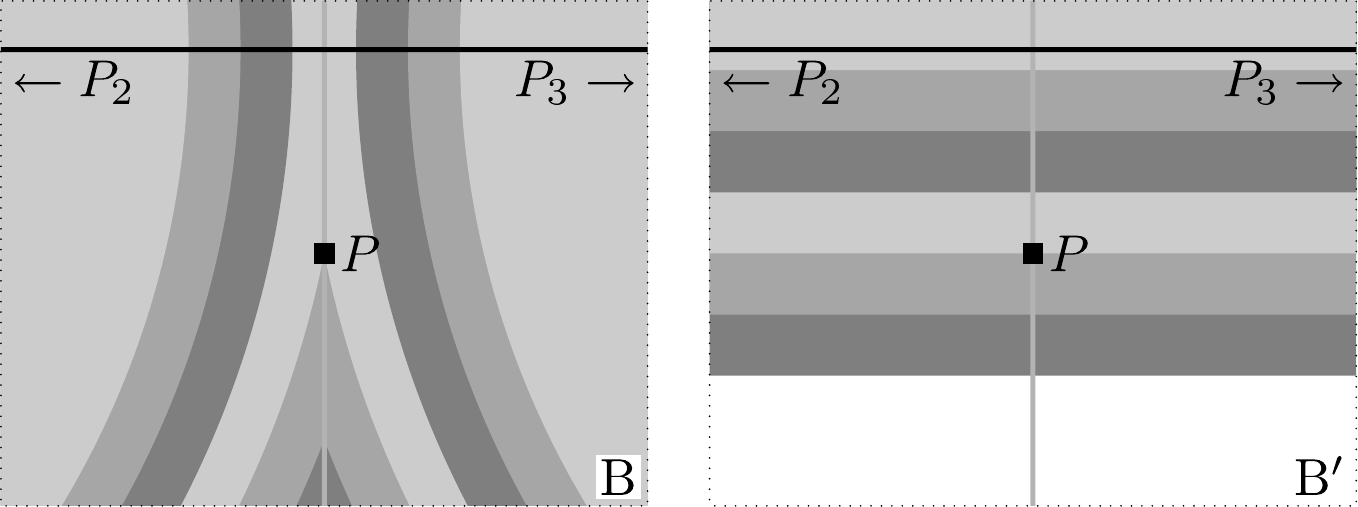}
\end{center}
\caption{Topologically regular situation around $P$ in picture B. The homeomorphic local change of variables makes the level set look as in picture B'.} \label{fig:fig2}
\end{figure}

Consider a point $P$ on the {\em interior of a Voronoi-edge}
$\Vor (P_i,P_j)$.
Suppose first that $P$ is different from the center of the segment $P_iP_j$. Then
there is no change in the topology of the lower level sets $\{ g \le g(P) \}$, since the set of level curves of $g$ is
topologically equivalent to a set of parallel lines. More
precisely: there exists a homeomorphism $\phi$ of a open
neighborhood of $P$ onto an open set in $\RR^2$ such that the
composed function $g \phi$ is a linear function. In this case we
call $P$ a {\em topologically regular point} of $g$, see figure \ref{fig:fig2}.
Suppose next that $P$ coincides with the center of the segment $P_2P_3$,
then it is possible to make a non-differentiable (but
homeomorphic) change of coordinates $\phi$, such that the composed
function  $g \phi$ is given by the formula: $ g(P) + x^2 - y^2$,
which defines a differentiable saddle point.
\begin{figure}[htb]
\begin{center}
\includegraphics[width=\textwidth]{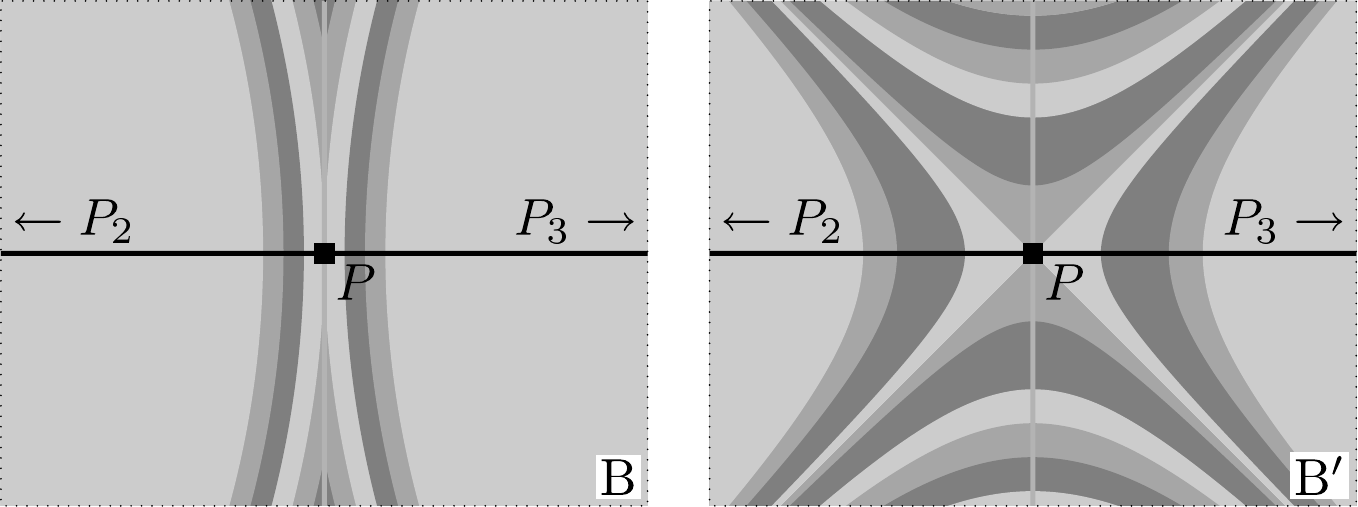}
\end{center}
\caption{At the point $c(P_2,P_3)$ of example B in figure \ref{fig:fig0} the function $g$ has a topological saddle point. The homeomorphic local change of variables makes the level set look as in picture B'.} \label{fig:fig3}
\end{figure}
In this case we call $P$ a {\em topological saddle point} of $g$
(figure \ref{fig:fig3}).
\par
After having considered neighborhoods of the edges of the Voronoi diagram, 
the remaining points to consider are the {\em vertices of the
Voronoi diagrams.} At a vertex $P$ several edges of the Voronoi diagram will meet.
We consider the following cases, related to the
behavior of the restriction of $g$ to the edges, containing  $P$.
\begin{itemize}
\item[i.] {\it$P$ is a maximum of this restriction.} We can now use a
non-differentiable (but homeomorphic) change of coordinates
$\phi$, such that the composed function  $g \circ\phi$ is given by the
formula: $ g(P) -x^2 - y^2 $, which defines a {\em local maximum}
of $g$ (figure \ref{fig:fig4}).
\begin{figure}[htbp]
\begin{center}
 \includegraphics[width=\textwidth]{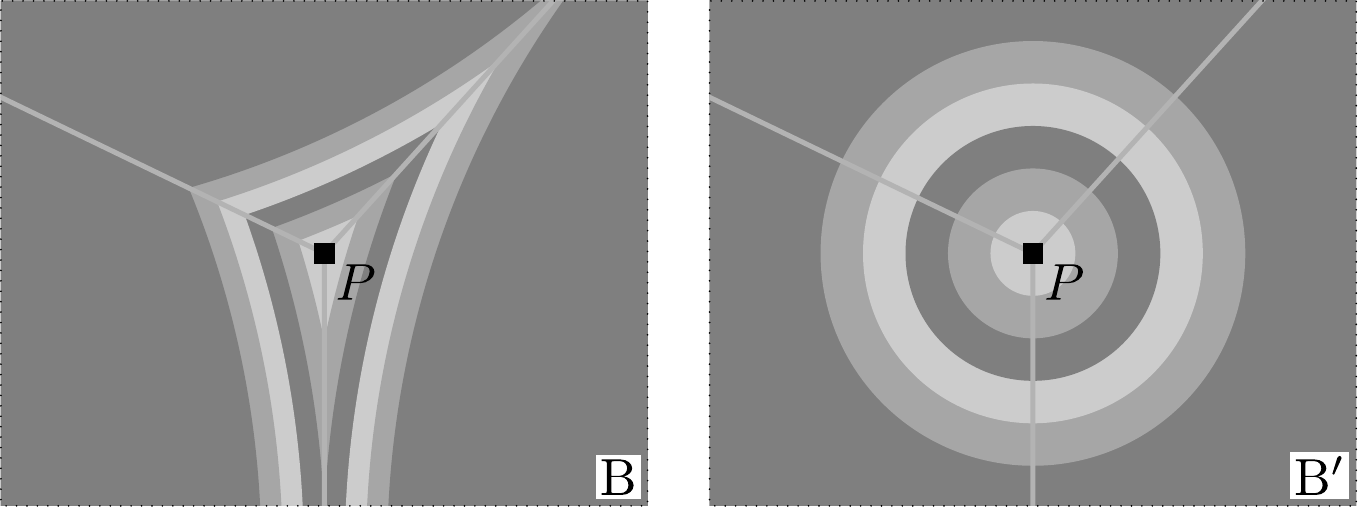}
\end{center}
\caption{A topological maximum} \label{fig:fig4}
\end{figure}
\item[ii.] {\it $P$ is a maximum on all but one of the adjacent edges.}
As figure \ref{fig:fig5} shows in this case there is a non-differentiable (but homeomorphic) 
change of coordinates that transforms $g$ into a linear function.
Again $P$ is a {\em topological regular point} of $g$.
\begin{figure}[h]
\begin{center}
  \includegraphics[width=\textwidth]{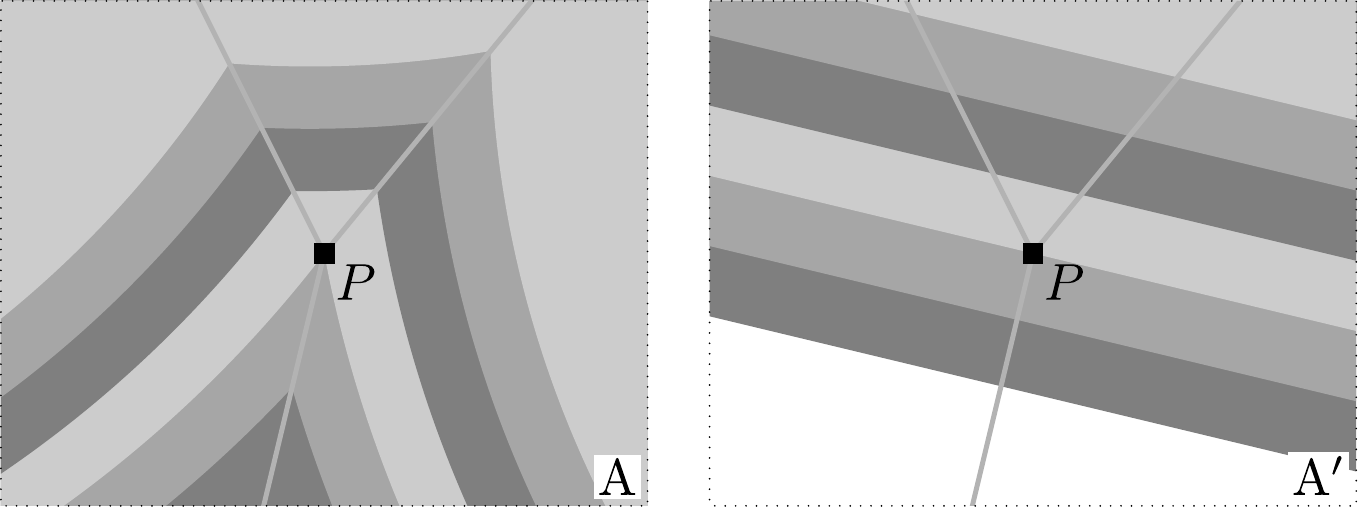}
\end{center}
\caption{A vertex, which is topologically regular} \label{fig:fig5}
\end{figure}
\end{itemize}
Note that cases other than (i) and (ii) do not exist. The reader may verify this
for three ajacent edges. That it holds in other cases ( 4 or more adjacent edges),
as well is proved in \cite{VoronoiDist}.
\par
The minima, the saddle points and the maxima of $g$ are called (topological)
critical points of index $0$, $1$, $2$ respectively. They occur as follows:
\begin{itemize}
\item[a.] minima in the centers of the Voronoi cells.
\item[b.] saddle points on the interior of edges of the Voronoi diagram,
\item[c.] maxima, which are vertices of the Voronoi diagram.
\end{itemize}
but as we have seen not any edge or vertex will give a critical point.
We will give a combinatorial criterium in \ref{DelMorse}.
The corresponding
values are called critical values. The other points are called
(topologically) regular. Besides the 
differentiable regular points in the interior of the Voronoi cells, we will meet
also topologically regular points, which lie on the interior of edges of the Voronoi diagram,
or are vertices of the Voronoi diagram.

\subsection{Morse formula}
The following Morse formula is a
non differentiable version
of the 'mountaineering equation'.
\begin{theorem}\label{sec:morse-formula}
Let $s_0$, $s_1$ and $s_2$ be respectively the number of (topological) minima,
saddle points and maxima of the  distance functions $d$ or $g$. We have:
\begin{equation*}
s_0 - s_1 + s_2 = 1
\end{equation*}
\end{theorem}
For the proof one can use the framework of Morse theory, which is well
known in differential topology. See Milnor \cite{Mi} or Hirsch
\cite{Hi}. In this article we use topological Morse theory. Morse theory was
carried over to the topological case by Morse himself
in the articles \cite{MR22:4071} and \cite{MR47:9631}. 
Voronoi diagrams are examples of stratified spaces.
Generalizations of Morse theory to stratified spaces by Goresky and MacPherson are discussed in \cite{GM}.
\par
For the proof of \ref{sec:morse-formula}
the main idea is to follow $\chi \{g \le t\}$ and to check the formula
\begin{equation*}
   s_0(t) - s_1(t) + s_2(t) =  \chi \{g \le  t\}
\end{equation*}
if $t$ increases from $0$ to infinity. One starts at $t=0$ with $N$
minima only. During the process there are a finite number of
values, where $\chi$ changes with $-1$ for a saddle point and $+1$
for a maximum. For $t$ large enough we have a contractible set,
which covers almost $\RR^2$ and $\chi =1$ in that case. For more
details we refer to \cite{VoronoiDist}.
\subsection{Delaunay tesselation and Morse poset} \label{DelMorse}
The Delaunay tesselation with respect to the point set $\colP1N$ is the
division of the convex hull $\CH(\colP1N)$ in polyhedra dual to the
Voronoi tesselation of the plane.
We now describe how to construct the Delaunay tesselation.
The vertices of the {\em Delaunay tesselation} $\Del(\colP1N)$ are the
points of $\colP1N$. There is an edge connecting two points of
$\colP1N$ if and only if their Voronoi cells share a common edge (it is not
enough that two Voronoi cells intersect only in one point).
The $2$-simplices are the triangles where $3$ (or more)
Voronoi cells come together.

\begin{figure}[htb]
\begin{center}
 \includegraphics[height=5cm]{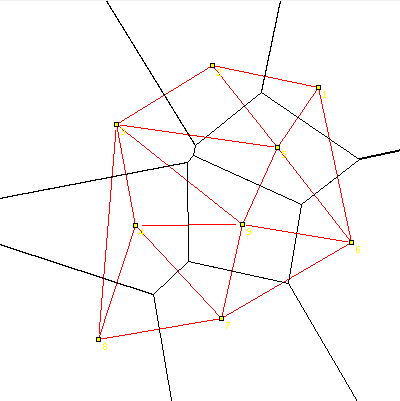}
\end{center}
\caption{Voronoi tesselation and and its dual Delaunay tesselation}
\label{fig:Delvor}
\end{figure}

\par

Note that if the convex hull of $\colP1N$ is $2$-dimensional and no four points lie
on a circle, then the Delaunay tesselation is in fact a triangulation of the convex hull
of $\colP1N$.
By abuse of language sometimes 'triangulation' is used, where in fact one has only a tesselation.
\par
{\bf Property.} The critical points of our distance function are exactly the intersection points between the Voronoi cells and the corresponding Delaunay cells. 

Saddle-points are the intersections
of two edges.  If the vertex 
of the Voronoi diagram is contained in the convex hull of the vertices of the corresponding Delaunay cell, then we have a maximum. 
For any point $P_i\in\colP1N$ the Voronoi and Delaunay cell intersect in the point itself. Indeed all the points $P_i$ are minima of the distance function. 
\par
So to each critical point there belongs a cell in the Delaunay tesselation, which we call
{\it active} or {\it critical}. These active cells constitute a poset: the Morse poset.
Note that the Morse poset is in general not a simplicial complex. It is only a poset. Recall that
a poset is a set with a partial order. In this case the partial order is inclusion. 
The Morse poset is a subset of the Delaunay tesselation.

\par
\begin{figure}[htpb]
\begin{center}
 \includegraphics[height=5cm]{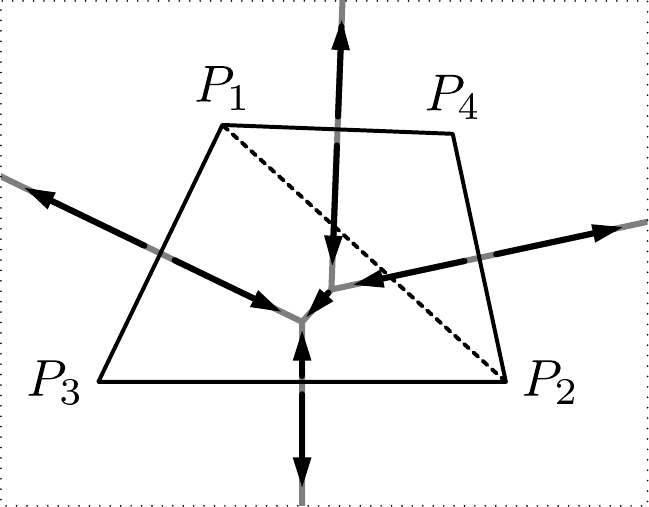}
\end{center}
\caption{Saddles-maxima graph.} \label{fig:fig17}
\end{figure}

It was shown in \cite{SieVanM} that the Morse poset is a  useful shape classifier.
Its 1-skeleton was called in \cite{VoronoiDist} the  {\em saddles-maxima graph}
of $\colP1N$, for short $\sm(\colP1N)$ (cf.\  figures \ref{fig:fig17} and \ref{fig:8}).
This saddles-maxima graph is known in computer science as the {\em Gabriel graph}, 
cf.\ \cite{GS} and  \cite{Ur}.
We leave it to the reader to show that in fact the saddles-maxima graph $\sm$
is connected. 

\subsection{Next to a discrete Morse function}\label{sec:next-discrete-morse}
To the Morse poset we can add extra
information about the topological cell attaching: just give a cell of the Delaunay tesselation the
value of the distance function $g$ at the corresponding critical
point. We call this a discrete function on the Morse poset.
\par
This is near to the concept of discrete Morse theory, which we
will describe in section \ref{sec-discr}. The starting point there is a
function on the simplices of a simplicial complex. 
In our case we could take the Delaunay tesselation and extend the discrete function
with values on the simplices which are not in the Morse poset.
\par
One can do this as follows  in the (above) two 3-points examples :
\begin{itemize}
\item Acute triangle. All cells belong to the Morse poset.
\item Obtuse triangle. The three vertices and two of the edges belong to the Morse poset;
assigne to the third edge and the 2-simplex the maximal value of $g$ on each of them (they turn out to be equal).
\end{itemize}
In more complicated  examples: follow the same strategy. Add a
missing $2$-cell together with a boundary $1$-cell as soon
as the two other boundaries are attached (give a value slightly
higher, but lower than the next critical value).
\par
We will treat this in more detail in section \ref{sec:morseposetToVf}. In theorem \ref{sec:from-morse-poset} we discuss the extension to a discrete Morse
function in all dimensions.
Discrete Morse functions are easy to implement on computers and in this case
they carry a lot of information about ``shape''.

\section{Power diagrams and tropical hypersurfaces.}
\label{sec-power}
In this section we show the equivalence of two definitions of
the power diagram. First we discuss the original notion of
Aurenhammer, by means of the power distance function, as used in computational geometry.
Then we show that these are equivalent to the
definitions by means of affine functions, as used in algebraic geometry. 
Note that we avoid genericity conditions in the treatment below.
\subsection{The definition of Aurenhammer}
Take a point set $\colP1N\subset\RR^n$. Throughout this article
we assume that $\dim(\CH(\colP1N))=n$.
Assign a weight $w_i$ to each point $P_i$.
Then write down the functions:
\begin{equation}\label{eq:20}
  g_i \colon \RR^n \rightarrow\RR\quad g_i(x)=\frac12\lVert x -P_i \rVert^2 - \frac12w_i\quad g(x) = \min_{1\leq i \leq N}g_i(x)
\end{equation}
Here we have used the following notation:
\begin{equation*}
  \lVert x\rVert^2 = \sum_{i=1}^nx_i^2
\end{equation*}
Note that $2g_i(x)$ is equal to the the square of the tangent from the point $x$ to the sphere around $P_i$ with radius
$r_i=\sqrt{w_i}$ (as soon as $x$ is outside that sphere). This value is in Euclidean geometry known as the power of a point with respect to that sphere. Points on the sphere have $2g_i(x)= 0$ and inside the sphere one has $2g_i(x) < 0$. 
\begin{figure}[ht]
\begin{center}
 \includegraphics[width=0.6\textwidth,viewport=0 300 400 580]{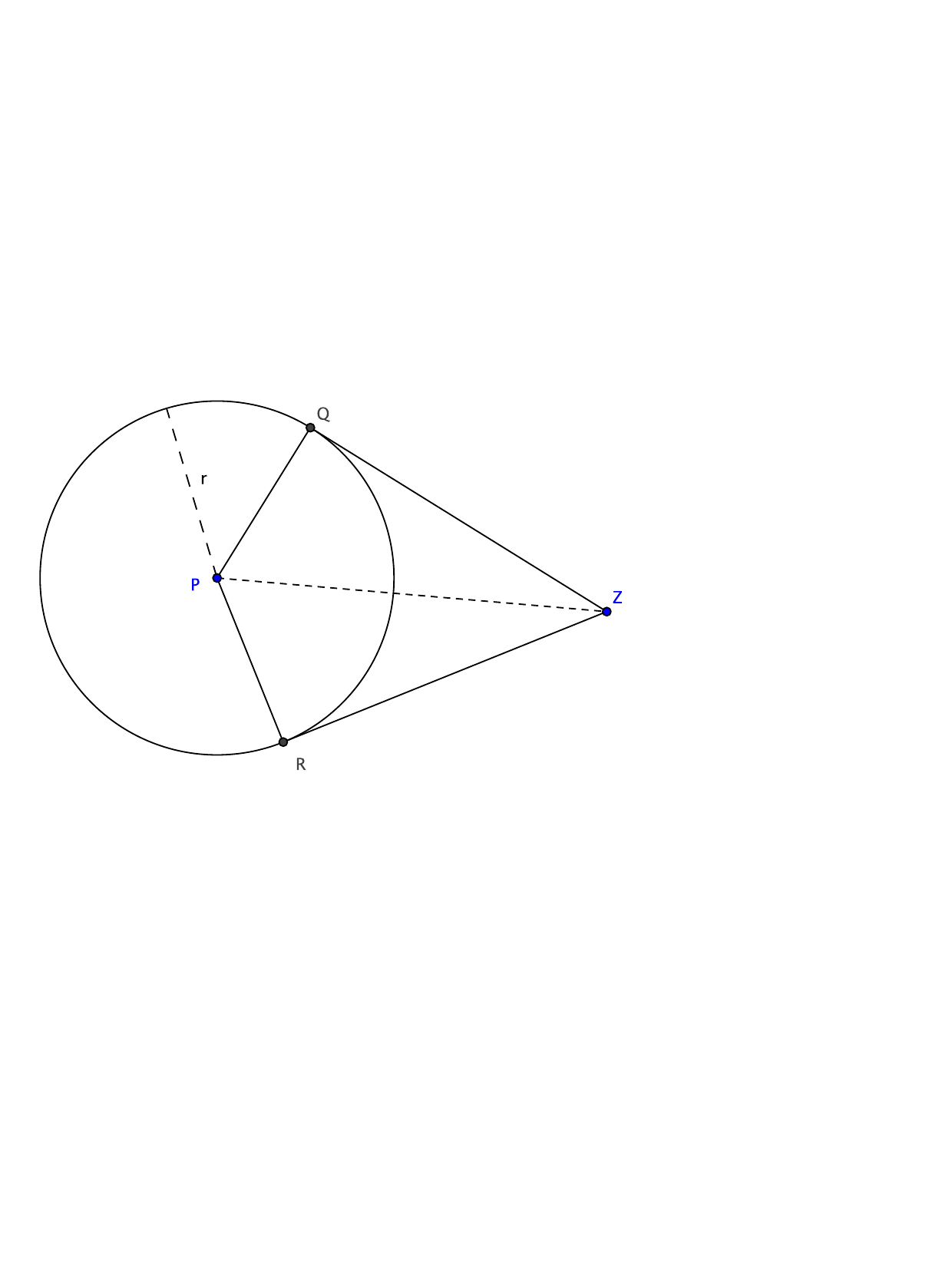}
\end{center}
\caption{Power of a point with respect to  a sphere} \label{fig:pointcircle}
\end{figure}

\begin{definition}
For a subset $\alpha\subset\colP1N$ the set $\Pow(\alpha)$ is the closure of
\begin{equation*}
  \{ x \in\RR^n\,\mid\,
   g_i(x)  = g(x) \,\,\, P_i\in\alpha\text{ and }
   g_j(x)  > g(x) \,\,\, P_j\not\in\alpha
\}
\end{equation*}
The sets $\Pow(\alpha)$ are called (power) cells.
\end{definition}
It is no restriction to assume that $w_i>0$. We may add some 
number to all of the functions $g_i$: the cells will not change. 
In fact we can add any function $h\colon\RR^n\rightarrow \RR$ to 
the $g_i$. We will get the same power cells.
\par
We repeat definition 1 in \cite{MR2040284}.

\begin{definition}
A \emph{polyhedral subdivision} ( resp.\ polytopal subdivision )
$\mathcal{T}$ of a polyhedron $K\subset\RR^n$
is a subdivision of $K$ in polyhedra ( resp.\ polytopes ) $K_i$, such that 
\begin{itemize}
\item The union of all sets $K_i\in\mathcal{T}$ is $K$.
\item If $K_i$ and $K_j$ are both in $\mathcal{T}$ then so is their
intersection.
\item Every compact subset $L$ of $K$ intersects only a finite number of
the $K_i$.
\end{itemize}
\end{definition}
\begin{theorem}
The sets $\Pow(\alpha)$ for $\alpha\subset\colP1N$ 
are a polyhedral subdivision of $\RR^n$. This polyhedral subdivision is
called the \emph{power diagram}.
\end{theorem}
\begin{proof}
The cells cover $\RR^n$. They are the intersection of a finite
number of half spaces of the form:
$\{ x\in\RR^n \,\mid\, g_i(x)\leq g_j(x) \}$. So they are polyhedra.
\end{proof}
Besides many similarities , there are at least two differences
between power diagrams and Voronoi diagrams.
\begin{itemize}
\item In a power diagram the cell of a point may well be empty.
It might happen that for some $i$ there is for
every $x\in\RR^n$ a $j=j(x,i)$ such that $g_j(x) < g_i(x)$.
\item Whereas in Voronoi diagrams, $P_i$ is always contained in its own cell,
in a power diagram the cell of $P_i$ might not be empty and still $P_i$
does not lie in its own cell.
\end{itemize}
We will see instances of these two phenomena in the examples below.
\par
\begin{figure}[htbp]
  \centering
  \includegraphics[width=0.4\textwidth,viewport=0 0 400 350]{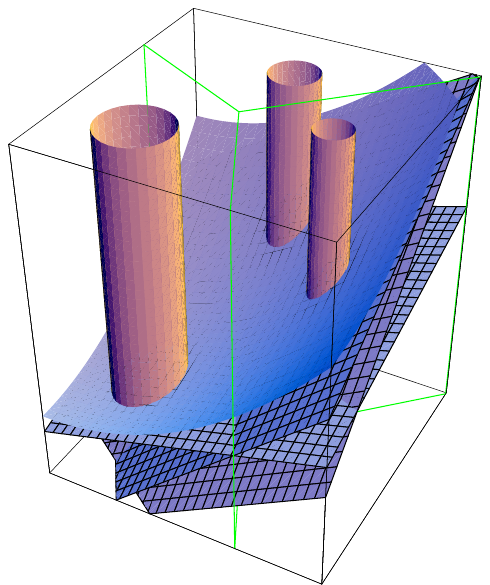}
  \includegraphics[width=0.4\textwidth,viewport=0 0 400 350]{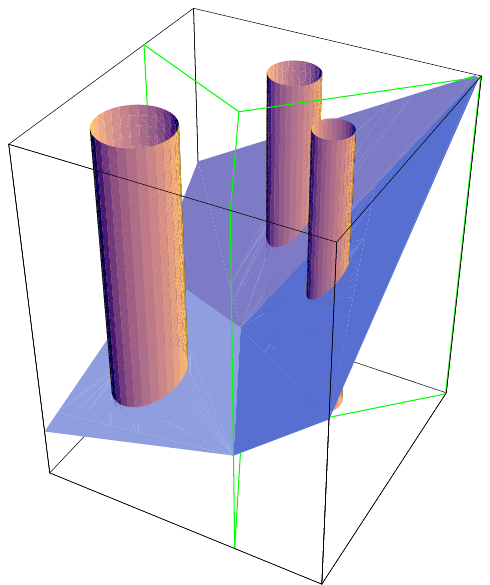}
  \caption{Power lifting to the parabola}
  \label{fig:12}
\end{figure}

Recall that the Voronoi diagram can be constructed using the upper convex
hull to the tangent planes to a parabola $x_0=\frac12\sum_{i=1}^n x_i^2 $ at the
points $\colP1N$. For power diagrams a similar construction exists.
Look at figure \ref{fig:12}.
We take cylinders with radius $r_i=\sqrt{w_i}$ around the lines $x=P_i$ in
$\RR^{n+1}$. The intersection of the cylinders with the parabola
lies in a hyperplane. We draw the hyperplanes in which they lie ( left
figure ) and we consider the upper convex hull of these hyperplanes (right figure).
The projection of its singular sides is the power diagram.
\par
The above construction is described in \cite{MR873251}, section 4.1.
When the $r_i=0$ the construction
reduces to the well-known construction of Voronoi diagrams. In that
case, of course, one has to take the tangent planes to the paraboloid.
\subsection{Digression on Tropical Geometry}

We will connect the theory of power diagrams with the concept of
tropical hypersurface in tropical geometry.
For this purpose we make a short digression on this subject. 
For surveys on the subject we refer to \cite{AG0306366} and \cite{Ga}.
\par
Tropical geometry is a  relatively new development in algebraic geometry that tries
to connect algebraic geometry problems with combinatorial questions on certain polytopes.
Recall that algebraic geometry studies varieties: the zero set of polynomials with 
real or complex coefficients in affine space.
\par
In tropical geometry one considers two new operations in $\RR$:
\begin{itemize}
\item[-]
tropical addition: $x \oplus y := \min(x,y)$, and
\item[-]
tropical multiplication: $ x \otimes y := x + y$.
\end{itemize}
With these two operations $\RR$ gets the structure of a topological semi-ring.
Such a tropical semi-ring is called a {\it min-plus algebra}.
\par
Polynomials in tropical geometry are defined in the usual way. 
The ``dictionary'' from algebraic geometry to tropical geometry works as follows:
The ordinary polynomial $$4x^3  + 4y^3  + 2xy + 7$$ has a tropical version:
$$\min \{ 3x+4 , 3y+4, x+y+2 , 7 \}.$$ 
\par
Tropical polynomials are piecewise linear concave functions on $\RR^n$
with integer coefficients. They come together with a vertex set defined by
the exponents; in the example $(3,0)$, $(0,3)$, $(1,1)$, $(0,0)$.
\par
\begin{figure}[htbp]
  \centering
  \includegraphics[width=0.6\textwidth]{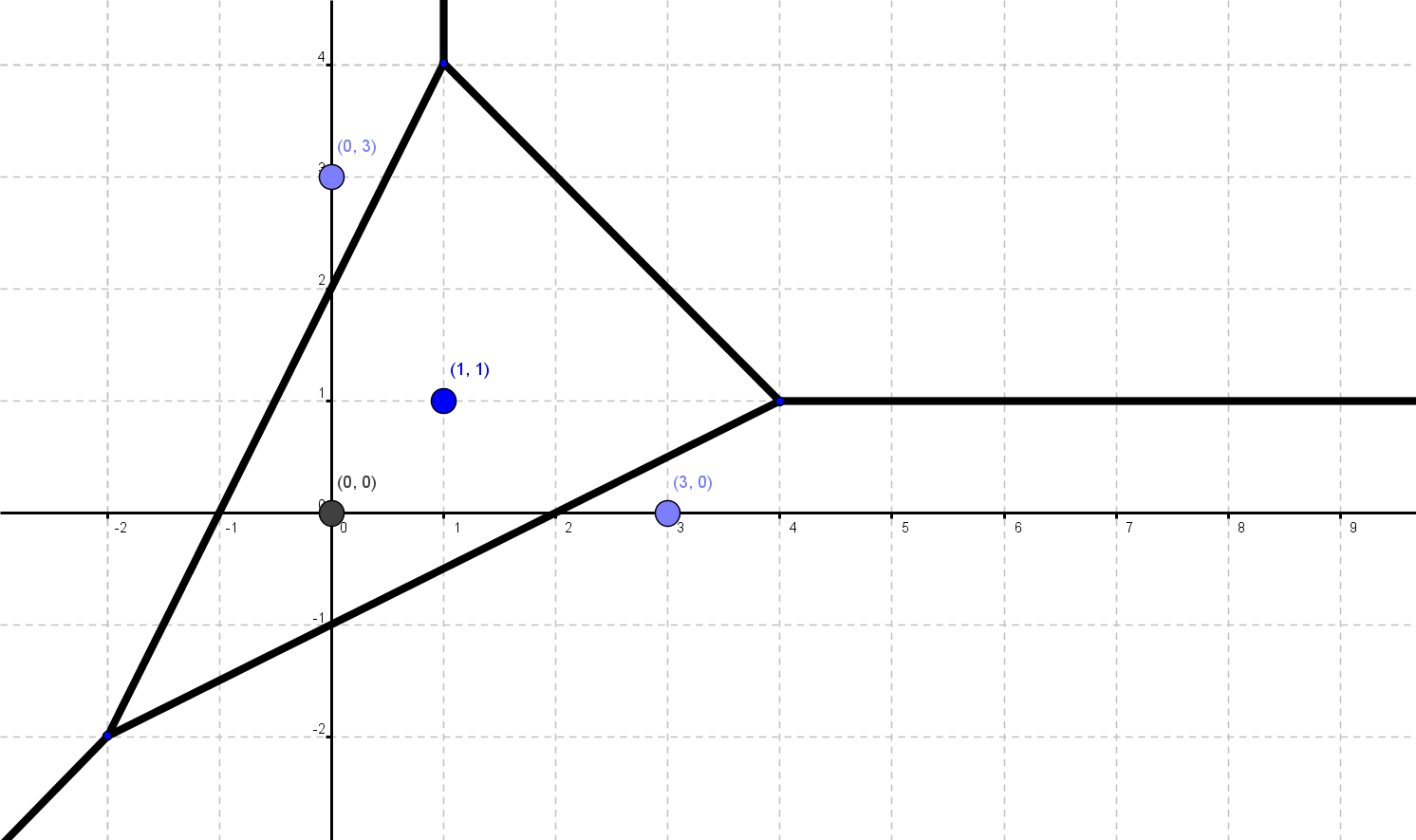}
 
  \caption{The corner locus for 
$\min \{ 3x+4 , 3y+4, x+y+2 , 7 \}$.}
  \label{fig:cornerloc4}
\end{figure}

The analogue of a variety in tropical geometry is the non-differentiability locus of 
the tropical polynomial, also called the corner locus of the concave function.
A corner locus 
 is drawn in figure \ref{fig:cornerloc4}. 
\par
 
The pictures one gets look similar to power diagrams.
Also tesselations of the polytope of the vertex set appear in a natural way
in tropical geometry.
\par
Because of the relation $\max(x,y) = - \min(-x,-y)$
one could also have used the operation $x \oplus y = \max(x,y)$ to define a tropical
semi-ring. Tropical hypersurfaces are in that context  piecewise linear convex functions. 
We will use this convention in the rest of this paper.
\par
The remarkable fact is that these tropical hypersurfaces also appear in a very different
way.
\par
Let $V \subset (\CC^\star)^n$  be an algebraic variety. Recall
that $\CC^\star = \CC - 0$ is the group of complex numbers under multiplication.
Let $\Log \colon (\CC^\star)^n \to \RR^n$ be the ``logarithmic moment-map'' defined by 
\begin{equation*}
\Log(z_1, . . . , z_n) = (\log \vert z_1\rvert , \cdots , \log \lvert z_n\rvert )\text{.}
\end{equation*}
Gelfand-Kapranov-Zelevinski defined the {\it amoeba} of an algebraic variety $V$
as the image  $A = \Log(V ) \in \RR^n$.
As one can see in the articles \cite{MR2040284}, \cite{Ga} and elsewhere that
the tropical hypersurfaces (drawn with the maximum convention) are ``spines'' of amoebas of algebraic varieties.
\begin{figure}[htbp]
  \centering
  \includegraphics[height=5cm]{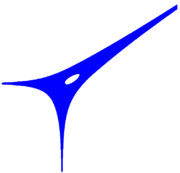}
 
  \caption{Amoeba of $3x^2 + 5xy + y^3 + 1$}
  \label{fig:amoeba}
\end{figure}
It is a deep result in tropical geometry established by Mikhalkin that the spine is 
a certain limit of the amoeba and carries all topological information
of the amoeba.
\par
For details we refer the interested reader to \cite{AG0306366}, \cite{Ga}, \cite{MR1264417} and \cite{Ga}.
\subsection{Power diagram as tropical hypersurface.}
As we hinted at in the previous section
the power diagram can also be defined using only affine functions.
Let as in equation (\ref{eq:20}):
\begin{equation*}
  g_i(x)=\frac12\lVert x -P_i \rVert^2 - \frac12w_i\quad g(x) = \min_{1\leq i \leq N}g_i(x)
\end{equation*}
Take the affine functions
$f_i(x)=\langle x , P_i \rangle + c_i$, where the coefficient $c_i$ is
\begin{equation}\label{eq:15}
  c_i = - \frac{\lVert P_i \rVert^2 -w_i}2 
\end{equation}
Consider the maximum of the $f_i$:
\begin{equation*}
  f(x)= \max_{i=1, \cdots , N} f_i(x)
\end{equation*}
Then $f$ and $g$ satisfy the following relations:
\begin{equation}\label{eq:2}
g_i(x) = \frac12\lVert x  \rVert^2  - f_i(x)  \quad   g(x) = \frac12\lVert x  \rVert^2  - f(x)
\end{equation}
The relations in equation (\ref{eq:2}) show that we can define the
power diagram solely using affine functions.
\begin{proposition}
For a subset $\alpha\subset\colP1N$ the set $\Pow(\alpha)$ is the closure of
\begin{equation*}
\{ x \in\RR^n\,\mid\,
f_i(x)  = f(x) \,\,\, P_i\in\alpha\text{ and }
f_j(x)  < f(x) \,\,\, P_j\not\in\alpha
\}
\end{equation*}
\end{proposition}
Using the affine definition it is immediate that the following two statements
that are false for Voronoi diagrams are true for power diagrams. These
two statements can also be found in \cite{MR873251}.
\begin{proposition}
The intersection of a hyperplane with a power diagram is again a
power diagram. The image of a power diagram under an affine map
$\RR^n\rightarrow\RR^n$ is again a power diagram.
\end{proposition}
So a power diagram can be sliced, after which we obtain a 
new power diagram. In fact, every power diagram in $\RR^n$ is a slice
of a Voronoi diagram in $\RR^{n+1}$.
\begin{proposition}
Let $\mathcal{T}$ be a power diagram in $\RR^n$. There is a Voronoi diagram
$\Upsilon$
in $\RR^{n+1}$ and a hyperplane $H$ such that $H\cap\Upsilon=\mathcal{T}$.
\end{proposition}
\begin{proof}
We will explicitly construct such an $\Upsilon$.
Let $\colP1N$ be the set of points in $\RR^n$ that determine $\mathcal{T}$.
Write in $\RR^{n+1}$:
\begin{equation*}
Q_i=(P_i,P_i^\prime ) \quad \bar{x}=(x,x^\prime)
\end{equation*}
We have 
\begin{equation*}
  \lVert \bar{x}-Q_i \rVert^2 = \lVert x - P_i \rVert^2 + (x^\prime - P_i^\prime)^2
\end{equation*}
So when $x^\prime=0$ we get
\begin{equation*}
  \lVert \bar{x}-Q_i \rVert^2 = \lVert x - P_i \rVert^2 +  P_i^{\prime 2}
\end{equation*}
Now let the power diagram be given by functions $g_i$ as in equation
\eqref{eq:20}. It is no restriction to assume that $w_i<0$, because
only the differences $w_i-w_j$ matter for the power diagram. 
Thus we can choose $P_i^\prime=\sqrt{-w_i}=r_i$. 
\end{proof}
\subsection{The Legendre transform}
In the 2-dimensional case we saw that
the Delaunay tesselation is the dual of the Voronoi diagram.
In a similar way there exists in any dimension a dual of the power diagram.
Using the tropical point of view
( with affine functions ) we can best explain the dual object using the
Legendre transform, see \cite{MR1301332}, which we explain below.
This dual object is also considered in \cite{MR1264417}, where it has
the name \emph{coherent triangulation}. 
\par
\begin{figure}[htbp]
  \centering
  \includegraphics[width=0.4\textwidth,viewport=0 0 400 320]{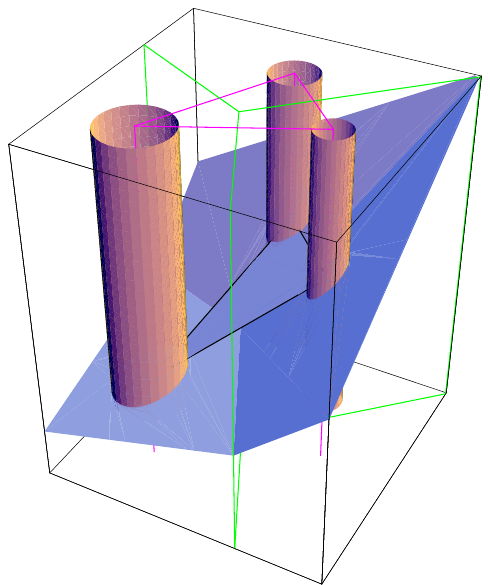}
  \includegraphics[width=0.4\textwidth,viewport=0 0 400 320]{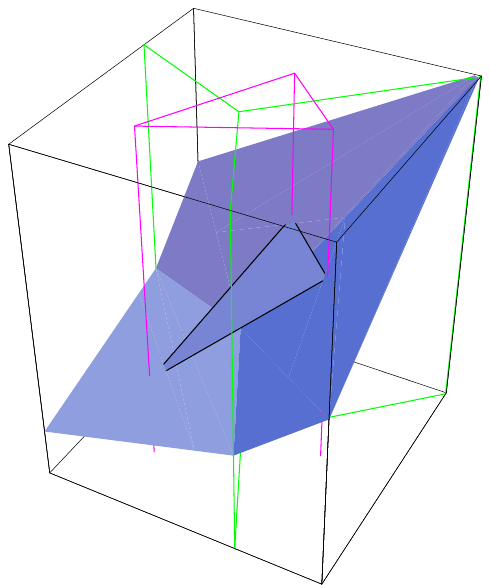}
  \caption{The power lift and the Legendre transform}
  \label{fig:11}
\end{figure}

We start with some definitions.
Let $P$ be a polyhedron in $\RR^{n+1}$. Let $v$
be a vector in $\RR^n$. The \emph{lower faces} of $P$
with respect to $v$ are those faces $F$ of $P$ such
that 
\begin{equation*}
  \forall x \in F\,\, \forall \lambda \in \RR_{>0} \colon x - \lambda v \notin P 
\end{equation*}
A polytopal subdivision of a polytope in $\RR^n$ is called
\emph{coherent} if it is the
projection of the lower faces of a polytope in $\RR^{n+1}$.
Not every polyhedral subdivision is coherent, see chapter 5 in \cite{MR1311028}, or
chapter 7 in \cite{MR1264417}. 
\begin{definition}
The Legendre transform of a convex function $f$, with domain
$D\subset \RR^n$ is
\begin{equation*}
\hat{f}(\xi)=\sup_{x\in D}\left(\langle \xi, x \rangle - f(x)\right)
\end{equation*}
When the supremum does not exists, we put $\hat{f}(\xi)=\infty$.
The domain $\Dom(\hat{f})$ of $\hat{f}$ are those $\xi$ for which
$\hat{f}(\xi) < \infty$.
\end{definition}
The Legendre transform of $f(x)=\frac12\lVert x \rVert^2$ is
the function itself. The Legendre transform of a linear
function $f_i=\langle x , P_i \rangle +c_i$ is $ < \infty$
only when $\xi =P_i$.
Theorem 2.2.7 of \cite{MR1301332} reads:
\begin{theorem}\label{sec:legendre-transform}
Let $f = \sup_{\alpha\in A}f_\alpha(x)$ be the maximum of a 
number of lower semi-continuous convex function. Then $f$ is
also a lower semi-continuous convex function. Furthermore
$\hat{f}$ is the infimum over all finite sums:
\begin{equation*}
  \hat{f}(\xi)=\inf_{\sum\lambda_\alpha\xi_\alpha=\xi,\lambda_\alpha>0,\sum_\alpha\lambda_\alpha=1}
\sum_\alpha\lambda_\alpha\hat{f_\alpha}(\xi_\alpha)
\end{equation*}
\end{theorem}
With theorem \ref{sec:legendre-transform} we calculate next
the Legendre transform of the function that determines the power
diagram. We have
\begin{equation*}
  f=\max_{1 \leq i \leq N} \langle x , P_i \rangle + c_i
\end{equation*}
The domain where $\hat{f} < \infty$ is the convex hull of 
the points $\colP1N$. We have
\begin{equation*}
  \hat{f_i}(P_i)=-c_i \Rightarrow \hat{f}(\xi)=\inf \left( -\sum_i\lambda_i c_i \right)
\end{equation*}
where the infimum is taken over all $\lambda_i$ such that 
\begin{equation*}
  \sum_{1\leq i \leq N}\lambda_iP_i = \xi\text{  and  }
  \forall i\colon \lambda_i \geq 0
\text{  and  } \sum_{i=1}^N \lambda_i = 1
\end{equation*}
It is no restriction to assume that not more than $n+1$ of the 
$\lambda_i$ are non-zero because any point in the convex hull
of the $P_i$ can be expressed as a sum of $\dim\CH(\colP1N)$ of 
the $P_i$.
The infima are best thought of in a geometric way. 
Take the convex hull $\CH( \{ ( P_i, -c_i) \}_{i=1}^N )$.
The lower faces form the
graph of $\hat{f}$ over the convex hull $\CH(\colP1N)$. We summarize
our discussion in a theorem, that can also be found in \cite{MR1264417}.
\begin{theorem}\label{thm:legendre-transform}
Let $\mathcal{T}$ be a power diagram. Let $\Upsilon$ be the
lower convex hull of the lifted points $(P_i, -c_i)$ wrt.\ to the vector
$( 0 , \cdots , 0 , 1 )$.
The polyhedral complex $\Upsilon$ is the graph of the Legendre
transform $\hat{f}$ of $f$.
The domain $\Dom(\hat{f})$ of $\hat{f}$ is
the convex hull $\CH(\colP1N)$.
\end{theorem}
The geometrical construction is directly related to the construction
of figure \ref{fig:12}. If we put the upper convex hull of theorem
\ref{thm:legendre-transform} in figure \ref{fig:12} we get figure
\ref{fig:11}.
What exactly happens here can best be understood by writing out the
equations. In figure \ref{fig:12} the graph of $x_0=\frac12\lVert x \rVert^2$
is drawn. Then the cylinders in the picture are
$\lVert x - P_i \rVert^2 = r_i^2$.
Or $\lVert x -P_i \rVert^2 = w_i$.
Hence the three planes of which we determine the upper hull are
\begin{equation}\label{eq:21}
  x_0=\langle x , P_j \rangle +\frac{r_j^2-\lVert P_j \rVert^2}{2}
\end{equation}
The triangle in figure \ref{fig:11} is
the plane $x_0=-\sum \lambda_i c_i$ where the $\lambda_i$ are
defined by $x=\sum\lambda_iP_i$.
In the point $P_i$ we have
that the graph of $f(x)=\max_if_(x)$ lies $\frac{w_i}2$ above the
paraboloid $x_0=\frac12\lVert x \rVert^2$, whilst the graph of
$\hat{f}(\xi)$ lies $\frac{w_i}2$ below the paraboloid there.
\par
Theorem \ref{thm:legendre-transform} gives a dual coherent triangulation
of the power diagram. 
\begin{definition}
Given a power diagram $\mathcal{T}$ specified by $N$ affine functions on $\RR^n$,
the coherent triangulation or (power) Delaunay tesselation of $\colP1N$ is the division of  $\CH(\colP1N)$ into
polytopes bounded by the set of points
$\xi\in\CH(\colP1N)$ where $\hat{f}$ is not differentiable.
\end{definition}
Since the coherent triangulation is the projection of the lower convex hull,
it is immediate that the sets bounded by the corner locus are in fact polytopes.
We will use `coherent triangulation' in relation with affine functions and `(power) Delaunay tesselation' in relation to the quadratic power distance. We also will omit `power'.
In the case when all the weights $w_i$ are zero this coherent triangulation
is of course the (classical) Delaunay tesselation.
\par
We next discus some general position arguments.
A number of points in $\colP1N$ in $\RR^n$ is affinely in general
position when the convex hull of each subset of length $n+1$ of
$\colP1N$ is full dimensional. For Voronoi diagrams  
it might still be that $n+2$ lie on a sphere. In
that case the Delaunay tesselation is not a simplicial complex. 
As condition for the (classical) Delaunay tesselation to really be a triangulation we need
that the points $\colQ1N$ are affinely general position in $\RR^{n+1}$ where
 $ Q_i = ( P_i , \frac12\lVert P_i \rVert^2 )$.
\begin{definition}
A power diagram $f=\max f_i$ is \emph{in general position}
when the points $Q_i$ are affinely in general position in
$\RR^{n+1}$. Here
\begin{equation*}
  Q_i=(P_i, -c_i)\text{, } 1 \leq i \leq N
\end{equation*}
\end{definition}
The following theorem is left to the reader.
\begin{theorem}
The Delaunay tesselation is a simplicial complex when the 
power diagram is in general position.
\end{theorem}
We can combine $f$ and $\hat{f}$ in the following function:
\begin{equation*}
  (\xi , x) \mapsto F(x,\xi )= f(x) + \hat{f}(\xi) - \langle x , \xi \rangle 
\end{equation*}
The function hides the \emph{Gateau differential} of $f$. Namely let
$h\colon\RR^n\rightarrow\RR$ be a function. If the limit
\begin{equation*}
  h^\prime(x ; \xi ) = \lim_{t \downarrow 0} \frac{h(x+t \xi ) - h(x)}{t}
\end{equation*}
exists, it is called the Gateau differential. See \cite{MR1301332}, theorem 2.1.22.
For a convex set $K\subset\RR^n$ the \emph{supporting function} of $K$ is
the function 
\begin{equation*}
  \xi \mapsto \sup_{x\in K} \langle x , \xi \rangle 
\end{equation*}
Theorem 2.2.11. in \cite{MR1301332} says that the Gateau differential 
$\xi \rightarrow f^\prime( x ; \xi )$ is the supporting function of the
set $\delta f(x)$ where 
\begin{equation*}
  \delta f(x)=\left\{ \mu \,\mid\, F(x,\mu) = 0 \right\}
\end{equation*}
Let us see what that means. Take $\xi$ to be one of the points of
$\colP1N$. Then $\hat{f}(P_i)= - c_i$. So the set of $x$ where
$F(x,\xi)=0$ are those for  which $f(x)=\langle x, P_i \rangle + c_i$.
\par
The Gateau differential is called Clarke's generalized derivative
in \cite{MR1481622}. It is stated in that article that
$\partial f(x)$ is the convex hull of the gradients of the
functions $f_i$ for which $f_i(x)=f(x)$. 
Let $\alpha$ be the set of points such that
$  f(x) = f_i(x)  \Leftrightarrow  P_i \in \alpha $.
So $x \in \Pow(\alpha)$.
Note moreover that these gradients are $\overrightarrow{xP_i}$.
So we get the convex hull of the points of $\alpha$, which is $\Del(\alpha)$.
\par
That last statement and the theorem 2.1.22 \cite{MR1301332} are all equivalent to
what is neatly formulated in proposition 1 of \cite{MR2040284}:
\begin{theorem}
There is a subdivision 
of the convex hull $\CH(\colP1N)$, dual to the power diagram $\mathcal{T}$.
The cell $\Del(\alpha)$ dual to $\Pow(\alpha)\in\mathcal{T}$ is
\begin{equation*}
\Del(\alpha) = \{ \xi \mid F(x,\xi) = 0 \quad\forall x \in\Pow(\alpha)\}
\end{equation*}
Reversely
\begin{equation*}
\Pow(\alpha) = \{ x \mid F(x,\xi) = 0 \quad\forall \xi\in\Del(\alpha)\}
\end{equation*}
\end{theorem}

\begin{remark}
One should make a difference between the Delaunay tesselation in the abstract sense, as subsets $\alpha$ of $\colP1N$,  and its geometric realization $\Del(\alpha)$. But in several cases, where we expect no confusion,
we will not distinguish between the two. So if $\alpha \subset \beta$ then we may say both `$\alpha$ is a face of $\beta$' or
 `$\Del(\alpha)$ is a face of $\Del(\beta)$'. Remark that (by duality) the power tesselation  is given by
 the same set of subsets of $\colP1N$, but with a different geometric realization by $\Pow(\alpha)$, which has complementary dimension.
 
\end{remark}

\section{The Morse poset}\label{morse}

\subsection{Introduction}
In section \ref{sec-2d} we studied in 2 dimensions the evolution of a set of wave fronts from a point set: growing discs with equal radius. We discussed critical points, the corresponding Morse poset, the changes in topology and the relation with Voronoi and Delaunay tesselations.
In this section we generalize to power distance functions. The corresponding Morse theory concerns the evolution of a set of growing balls in $n$-space with different radii. The union of fixed set of balls occur in models of molecules and are studied by e.g.~Edelsbrunner \cite{UnionOfBalls}. A related approach via alpha-shapes in $\RR^3$ can be found in \cite{DGJ}, \cite{GJ1} and \cite{GJ2}  . Further related developments are in \cite{Chazal}. 
  
\par
Let us  remark here that the topology of level sets of distance functions has been studied before  with
classical results that are not restricted to a finite set of points. In Riemannian geometry see for example \cite{Grove} and \cite{Cheeger};
in non-smooth analysis \cite{Clarke} and more recently in shape reconstruction \cite{Chazal}. We will use also the theory of continuous selections \cite{MR1481622}.

\par
We will define next Morse critical points in context of continuous functions. For topological Morse theory we refer to \cite{MR22:4071}.
\begin{definition}
Two functions $f_1 \colon U_1 \to \RR$  and $ f_2 \colon U_2 \to \RR$  are called \emph{topologically equivalent} if there exist a homeomorphism $\phi \colon  U_1 \to U_2$ such that $f_1 = f_2\circ\phi$.
A function  $f\colon M \to \RR$ is \emph{topologically regular} at $x \in M$ if there
exists a neighborhood $U$ of $x$ such that the restriction $f_{\mid U}$ is topologically
equivalent to a linear function.
A point $x$ is called a \emph{critical point} of $f$ if $f$ is not topologically regular at $x$.
The critical point is called \emph{topologically  isolated} if there exists a neighborhood
$U$ of $x$ such that all points $y \in U$,$y \ne x$ are regular points of $f$.
A topologically isolated critical point $x$ of  $f\colon M \to \RR$  is a
\emph{topological Morse point of index $d$} if there exists a neighborhood $U$ of $x$ such that the restriction $f|U$ has the property that for $\epsilon$ small enough the set
$\{x \in M \,\mid\, f(x) \le x+ \epsilon \}$
 is homotopy equivalent to  
$\{x \in M |f(x) \le x- \epsilon \}\, \cup\,  d-\text{cell}$. 
$f$ is called a topological Morse function if all critical points are topological Morse.
\end{definition}
\begin{example}
Standard examples $f\colon\RR^n \to \RR$ are
\begin{itemize}
\item classical Morse point; \cite{Mi}: \\
$ f(x) =   -(x_1^2 + \cdots + x_d^2) + (x_{d+1}^2 + \cdots + x_n^2)$
\item Morse point of continuous selections; \cite{MR1481622}:\\
$ f(x) =  \min( x_1,\cdots,x_d, -( x_1 + \cdots + x_d) )  + (x_{d+1}^2 + \cdots + x_n^2)$
\end{itemize}
\end{example}
We will omit the adjective ``topologically'' in the sequel and simply write: regular point,
Morse point, and so on.
\subsection{Critical points of the power distance function $g$}
Consider again the function 
\begin{equation*}
  g(x)=\min_{ 1 \leq i \leq N} g_i(x)
\end{equation*}
Critical points $x$ of $g$ are related to the generalized derivative $\partial g(x)$. 
For this concept and related notions we refer to \cite{Clarke}. 
Let $x$ be a point of the power diagram; then there is a set $\alpha$ such that
\begin{equation}\label{eq:powalph}
 g(x) = g_i(x)  \Leftrightarrow  P_i \in \alpha 
\end{equation}
So $x \in \Pow(\alpha)$.
Note that $\grad g_i(x) = x - P_i = \overrightarrow{xP_i}$.
The generalized derivative $\partial g(x)$  is the convex hull in the vector space $\RR^n$ of
the end points of the gradient vectors, so $\partial g(x)= \CH(\overrightarrow{xP_i})$.
By an affine transformation (sending the origin to $x$) one can identify
$\partial g(x)$ with $\CH(\alpha)= \Del(\alpha) $.
\par
Analogously to what was done in \ref{subs-crit} there are now three cases to consider:
\begin{enumerate}
\item[\bf Inside] $0\in \Interior(\partial g(x))$, hence: $x$ lies in $\Interior(\CH(\alpha))$,
or
\item[\bf Outside] $0\not\in \partial g(x)$, hence: $x$ lies outside $\CH(\alpha)$, or
\item[\bf Border] $0$ lies on the relative border of $\partial g(x)$, hence $x$ lies on the
border of $\CH(\alpha)$.
\end{enumerate}
\begin{figure}[htb]
\begin{center}
\includegraphics[width=\textwidth]{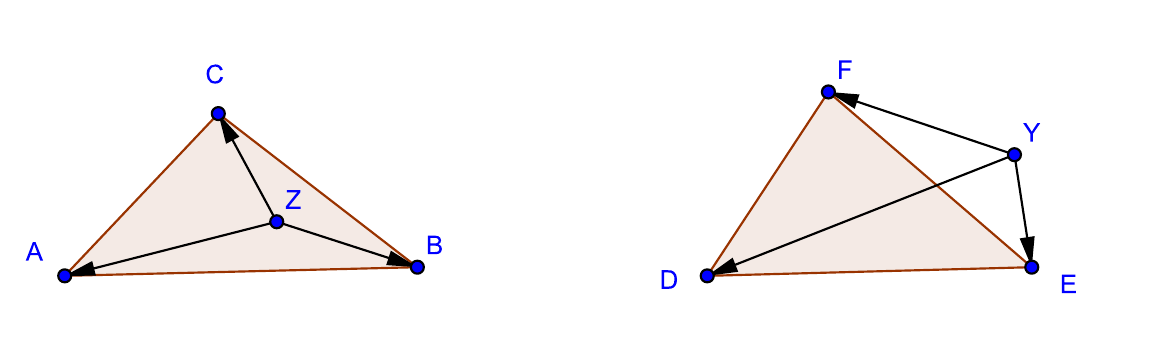}
\end{center}
\caption{The cases ``Inside'' and ``Outside''}
\label{fig:critical}
\end{figure}

``Inside'' implies that $g$ has a critical point by the standard theory of continuous selections 
theorem 2.2 in \cite{MR1481622}.
Note that this means that $x$ must be the (unique) intersection point of
$\Pow(\alpha) \cap \Del(\alpha)$;
both sets lie in orthogonal linear subspaces of complementary dimension.
Again compare this with the situation in section \ref{subs-crit}.

\par
Using the techniques of \cite{JS} 
one can see that 
``Outside'' implies that $g$ is regular at $x$. 
\par
The case ``Border'', which is not covered by the general theory, we will
treat it below. It should be noted that no general statements can be made 
about this case, but distance functions are special, see also \cite{Grove}, and for them some conclusions can be drawn.

Consider the Delaunay tesselation (dual to the power diagram).
In the special case of the Euclidean distance function we have defined
in section \ref{sec-2d}
a poset of active simplices in the Delaunay tesselation. They were defined
using the critical points of $\min g_i$. We repeat the procedure for
the power distance function.
\par
\begin{figure}[htbp]
  \centering
  \includegraphics[width=7cm,viewport=0 0 400 320]{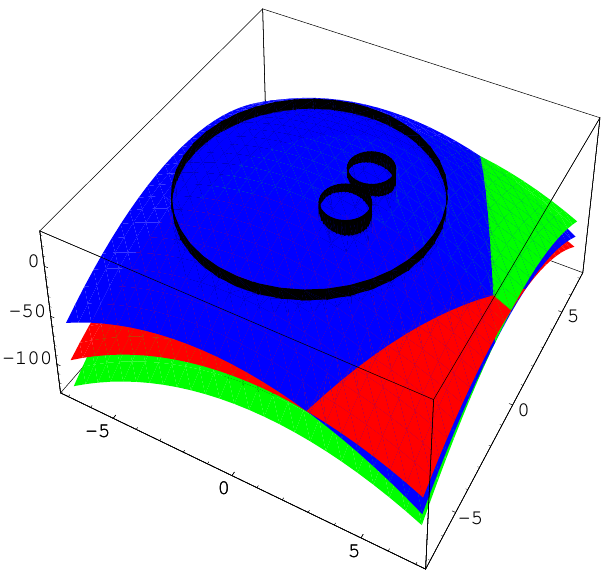}
  \caption{Typical counterintuitive example}
  \label{fig:13}
\end{figure}

Also, for a polyhedron $P$ in $\RR^{n+1}$ the relative interior is the interior
of the polyhedron as a subset of the affine span $\Aff(P)$ of $P$.
If $\dim\CH(\alpha)> 0$ then put
\begin{equation}\label{eq:24}
  c(\alpha)=\Sep(\alpha)\cap \Aff( \alpha )
\end{equation}
where $\Sep(\alpha)$ is defined as
\begin{equation*}
\Sep(\alpha)= \{ x\in\RR^n\,\mid\, \forall i,j\in\alpha \quad
\forall k\ni\alpha{}:f_i(x)=f_j(x)\not=f_k(x) \}
\end{equation*}
For a vertex put $c(\{ i \}) = P_i$. 
Recall that for a convex set $A$ a point $x$ is contained in the relative interior of 
$A$ if $x$ is an interior point of $A$ relative to the affine span $\Aff(A)$ of $A$.
\begin{definition}\label{def:active}
We call the cell in the Delaunay tesselation \emph{active} or \emph{critical} if
\begin{itemize}
\item $c(\alpha)$ is contained in the relative interior
$\Interior( \CH( \alpha ) )$ of $\CH(\alpha)$, and
\item $g( c(\alpha) ) = g_i(c(\alpha))$ for all $P_i \in \alpha$ and
 $g( c(\alpha) ) < g_j(c(\alpha))$ for all $P_j \not \in \alpha$.
\end{itemize}
\end{definition}
With this definition it might well happen that
 (in contrast with the Voronoi case) a vertex is not
active. Here is a typical counter intuitive example.
Take three points, say $P_1=(-1,0)$, $P_2=(1,1)$ and $P_3=(1,-1)$.
Then place a large circle around $P_1$: $r_1=5$. Place two smaller
circles around $P_2$ and $P_3$: $r_2=\sqrt\frac34$ and $r_3=\sqrt\frac56$. 
If we look from below at the graph of $\min g_i$ we see a picture
alike the one in figure~\ref{fig:13}.
\vspace{-1cm}
\begin{definition}
The \emph{Morse poset} of 
$\{ ( P_i, w_i ) \}_{i=1 , \cdots , N}$ consists
of the active subsets of the Delaunay tesselation. 
\end{definition} 
\begin{figure}[ht]
\begin{center}
\vspace{2cm}
 \includegraphics[width=0.3\textwidth,viewport=0 150 400 420]{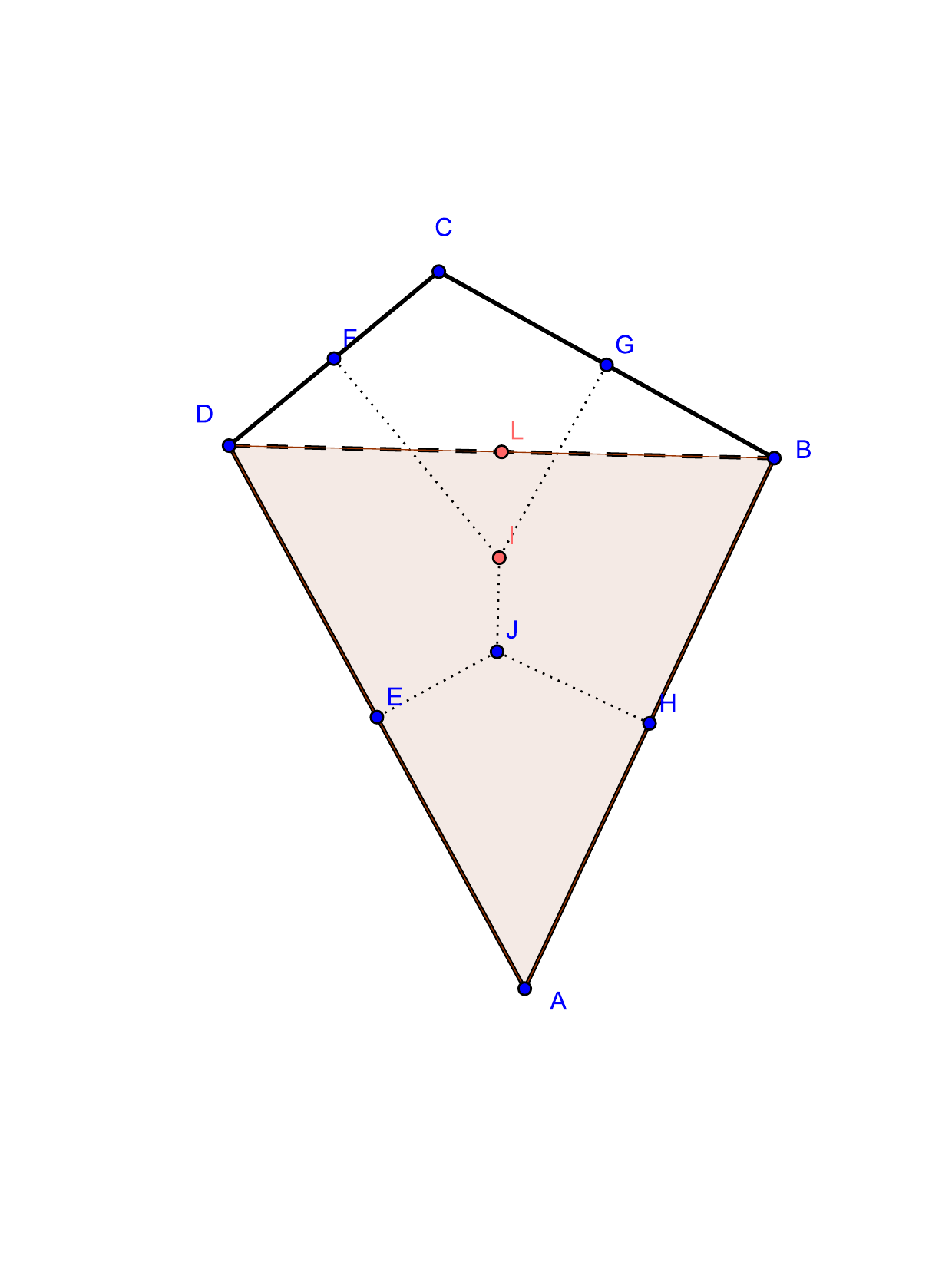}

\end{center}
\caption{Morse poset, consisiting of four $0$-cells, four $1$-cells and one $2$-cell.} \label{fig:4hoek}
\end{figure}
\par
Below we'll use the notation:
\begin{equation*}
g_\alpha(x) = \min_{i\in\alpha}g_i(x)\text{ and }
f_\alpha(x) = \max_{i\in\alpha}f_i(x)
\end{equation*}
\begin{theorem}
There is a one to one correspondence between critical points of $\min_i g_i$ 
and active cells in the Delaunay tesselation. An active cell of dimension
$d$ corresponds to a critical point of index $d$.
\end{theorem}
\begin{proof}
We will first rule out the case where $d=0$. Because the functions $g_i$ are all convex
with minima at their respective vertices,
a non-degenerate minimum can only occur on a vertex and the vertex will be active.
Conversely, if a vertex $P_i$ is active, we will have $g_j(P_i)>g_i(P_i)$ for all
$j\not= i$. And hence $g$ will have a non-degenerate minimum there.
\par
From now on assume that $d>0$.
\par
Let us again have $x$ as in equation \ref{eq:powalph}.
Note that equation \ref{eq:powalph} implies an open condition: in a neighborhood $x$ in $\RR^n$
the function $g$ is entirely determined by the functions $g_i$, with $i\in\alpha$. Rephrasing
that: there is a neighborhood $U$ of $x$ in $\RR^n$ such that $g_{\mid U}=g_{\alpha\mid U}$.
\par
As a consequence, $g$ will have a critical point of index $d$ at $x$ if and only if
$g_\alpha$ has a critical point of index $d$ at $x$.
We consider again the three cases ``Inside'', ``Outside'' and ``Border'' defined above.
``Inside'' implies that $x=c(\alpha)$, and hence $\alpha$ is active and
 $g$ has a critical point of index $d=\dim(\alpha)$ at $c(\alpha)$.
In more detail: choose coordinates $y=(x_1,\cdots,x_d)$ and $z=(x_{d+1},\cdots,x_n)$,
such that $\Aff(\alpha)$ is defined by $ x_1=\cdots =x_ d = 0$ and
$\Sep(\alpha)$ is defined by $ x_{d+1}=\cdots =x_ n = 0$.
Next: 

$$
g_{\alpha}(x) = \min_{i \in \alpha}{g_i(x)}
 = \min_{i \in \alpha}(1/2 \lVert y - P_i\rVert^2 - \frac12w_i +
\frac12 \lVert z\rVert^2 ) = $$  $$
= \min_{i \in \alpha}(1/2 \lVert y - P_i\rVert^2 - \frac12w_i) +
 \frac12 \lVert z\rVert^2
$$
So we have separation of variables. The first term defines a maximum in
the origin of $\RR^d$ and is equivalent to $-\frac12\rVert y\rVert^2 + C$.
So $g$ is equivalent to the standard Morse function of index $d$.
\par
As we remarked above``Outside'' implies that $g$ is regular at $x$.
Let us concentrate on the case ``Border''. We have seen that $\partial g(x)$ is a
polyhedron spanned up by the vertices in $\alpha$. 
\par
Let $\beta$ be the smallest face of $\alpha$ such that $x\in\CH(\beta)$.
We have that $x\in\CH(\beta)\subset\Aff(\beta)$, and $x\in\Pow(\alpha)\subset\Pow(\beta)$.
It follows that $x=c(\beta)=c(\alpha)$. We thus have to look at the level
sets of $g_\alpha$ near $x=c(\beta)=c(\beta)$.
\par
We claim that they are topologically
equivalent (i.e. homeomorphic) to the level sets of a linear function.
\par
In fact we can restrict attention to the behavior of $g_\alpha$ on $\Aff(\alpha)$. 
On the directions transversal to $\Aff(\alpha)$ and passing through
$c(\alpha)$ $g_\alpha$ is just a parabola with a minimum at $c(\alpha)$.
\par
But when we restrict to $\Aff(\alpha)$ the only singularity $g_\alpha$ can have at $c(\alpha)$ 
is a maximum because $g_\alpha$ attains its maximal value on $\Del(\alpha)$ at $c(\alpha)$.
We thus have to prove that $g_{\alpha\mid\Aff(\alpha)}$
does not have a local maximum at $c(\alpha)$.
\par
On $\Pow(\beta)$, whose only point of intersection with $\Pow(\alpha)$ is $c(\alpha)$,
we will have that $g_{\alpha\mid\Aff(\alpha)}\geq g(c(\alpha))$. In fact on $\Pow(\beta)$ the
value of $g_\alpha$ will ``wander off to infinity''.
Hence $g_{\alpha\mid\Aff(\alpha)}$ cannot have a local maximum at $c(\alpha)$.
\end{proof}

\subsection{Morse formula}
We have the following Morse formula for the power distance function.
\begin{theorem}
Let $s_i$ be the number of critical points of index $i$ of $g$. We have:
\begin{displaymath}
  \sum (-1)^i s_i = 1
\end{displaymath}
\end{theorem}
\begin{proof}
$g$ is a topological Morse function. In that case,
as $t$ grows, $g$ passes through a number of non-degenerate critical values.
When $g$ passes a critical value of index $i$, an $i$-cell gets attached \cite{Mi}.
In between we apply the (topological) regular interval theorem.
For each intermediate function value $t$ we have therefore :
\begin{equation*}
 \chi( \{g(x) \le t \}) = \sum (-1)^i s_i(t) 
\end{equation*}
\end{proof}
Finally we remark that though power diagrams are affinely defined,
activity is not retained under affine volume preserving linear
transformations. 
Through such transformations the power diagram and
its dual Delaunay tesselation do not change. The Morse poset though,
can change drastically after an affine transformation.
\section{Discrete Morse theory}\label{sec-discr}
\subsection{Introduction.}
As seen before the Morse poset is a subset of the Delaunay tesselation.
There is a natural function on this poset, which describes the order of cell attaching:
give each active cell the value of $g$ in the critical point.
An obvious question is how to extend this function to the full Delaunay tesselation.
In this way one enters the framework of Forman's discrete Morse theory.
A very good introduction to that subject is \cite{MR900443}. 
Descriptions of the theory exist in several categories, the most general
is CW-complexes, cf.\ \cite{MR1612391}. More common and more accessible for an outsider is the description for  simplicial complexes. We prefer to work inside the polyhedral category and avoid the CW-description.
\par
Forman considers a function $h$ on all cells of a simplicial complex, satisfying certain properties (see (\ref{eq:22}) and (\ref{eq:23}) below).
He also considers the concept of critical or active cell.
Following the evolution of values of $h$ gives a protocol for constructing
the simplicial complex from a set of points. This works in such a way, that non-critical
cells are attached in pairs and do not change the topology, while critical cells are attached
separately and therefore change the topology.
\par
In section \ref{sec:next-discrete-morse} we have described how to proceed from the distance function $g$ to a discrete Morse function that describes topology and shape in a convenient way.
In general there are two problems for the extension of our distance function $g$ to a
discrete Morse function $h$:
\begin{itemize}
\item The Delaunay tesselation is not always a simplicial complex. This already happens in the plane.
\item Some geometric defined candidates for $h$ do not always satisfy the Forman conditions. This is related to the fact that, unlike Morse functions on a manifold, Forman discrete Morse functions are not dense in the set of all functions on a simplicial complex. 
\end{itemize}
We will deal with both these issues and end up in section \ref{sec:morseposetToVf} with a discrete Morse function $h$, which extends $g$.
\subsection{Morse functions on Delaunay tesselation.}\label{sec:morse-funct-dela}
We first give the definition of a discrete Morse function and show next that the Delaunay tesselation has the property that for every two faces $\alpha\subset\gamma$ with 
$2+\dim(\alpha)=\dim(\gamma)$ there are two faces $\beta$
of the Delaunay tesselation in between $\alpha$ and $\gamma$.
This property will be important for constructing the discrete Morse function.
\par
Let $\mathcal{T}$ be a polyhedral subdivision. 
\begin{definition}
A function $h\colon\mathcal{T}\rightarrow\RR$ is called a \emph{Forman discrete
Morse function} if for all $\beta\in\mathcal{T}$
\begin{equation}\label{eq:22}
\# \{ \alpha\in\mathcal{T}\,\mid\, 1 + \dim(\alpha)=\dim(\beta)\,\,\,
 \alpha\subset\beta
 \,\,\, h(\alpha)\geq h(\beta) \} \leq 1
\end{equation}
and 
\begin{equation}\label{eq:23}
\# \{ \alpha\in\mathcal{T}\,\mid\, \dim(\alpha)=1+ \dim(\beta)\,\,\,
 \beta\subset\alpha
 \,\,\, h(\alpha)\leq h(\beta) \} \leq 1
\end{equation}
In case both numbers are zero for some $\beta\in\mathcal{T}$, $\beta$ is called
\emph{critical}. 
\end{definition}

\begin{example}
In figure \ref{fig:discMorse} we show a discrete Morse function. The values are indicated by  numbers. Since they are choosen all different, we also can use them to indicate the faces. The critical faces are:
  \begin{itemize}
  \item 0-dimensional: 1,2,4,6,10,11;
  \item 1-dimensional: 3,5,14,15,16,17;
  \item 2-dimensional: 18
  \end{itemize}
  Non critical faces occur in the following pairs:
  (8,9) and (12,13).\\
  If one interchanges the labels 3 and 1 then we don't satisfy the definition of discrete Morse function.

\begin{figure}[htbp]
  \centering
  \includegraphics[width=0.5\textwidth,viewport= 100 200 400 660 ]{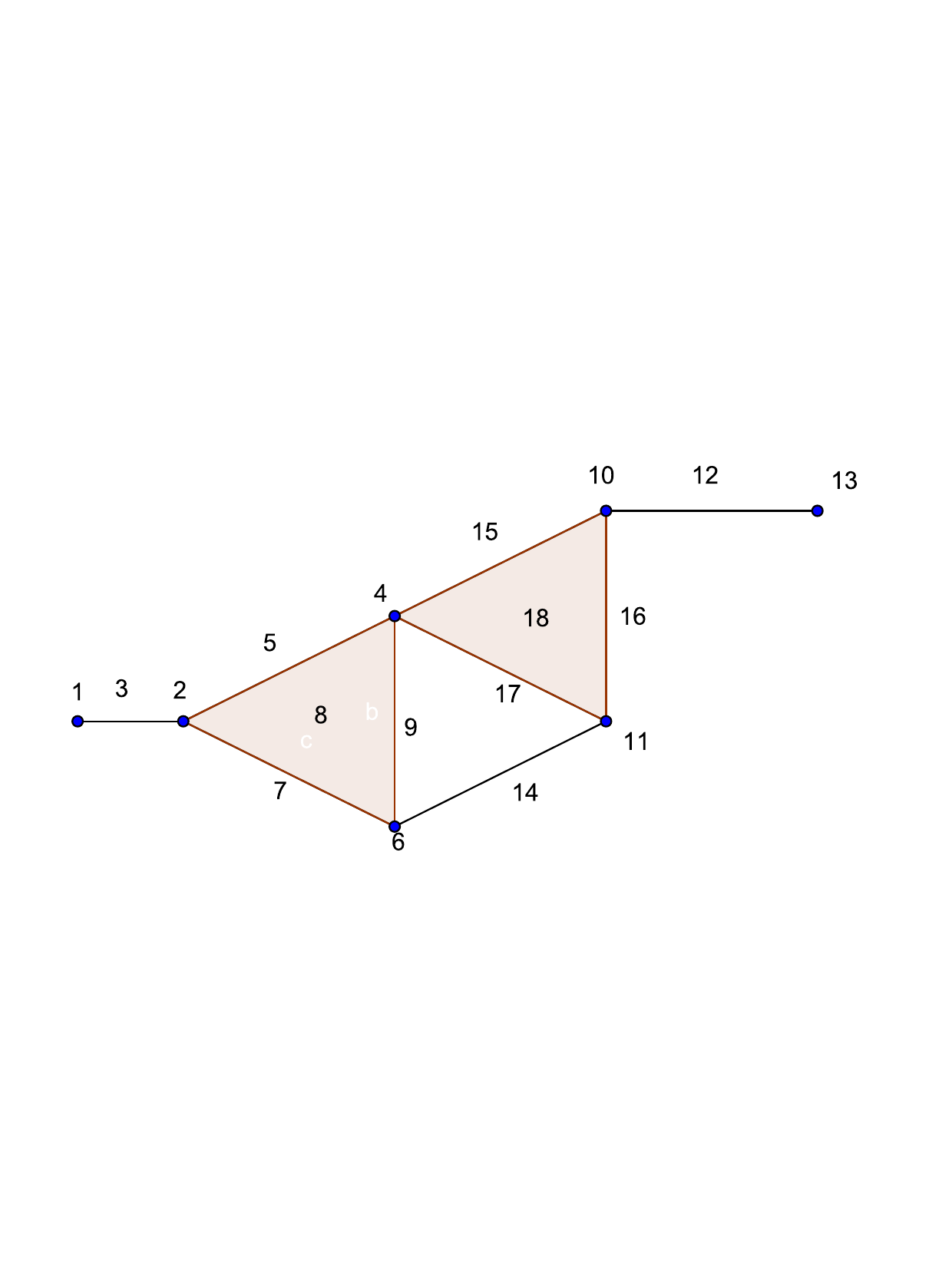}   
  \caption{A discrete Morse function}

  \label{fig:discMorse}
\end{figure}
\end{example}
We define a Morse poset of a discrete Morse function $h$ as the set
of critical faces of $h$. 
\par
If the power diagram $\mathcal{T}$ is in general position  then we know that it is a simplicial
complex.
Discrete Morse function on a simplicial complex have a very nice
property, see lemma 2.6 in \cite{MR1612391}.
If $\alpha$ is a face of $\beta$ and $\beta$ is a face
of $\gamma$. Label the vertices as follows:
\begin{equation*}
\alpha=\colP1{k}\quad\beta=\colP1{k+1}\quad\gamma=\colP1{k+2}
\end{equation*}
Because we are dealing with a simplicial complex there is another face
between $\alpha$ and $\gamma$:
$\delta=\{ P_1, \cdots ,P_k, P_{k+2} \}$:
\begin{equation*}
\UseTips
\xymatrix @C=1em @R=1em {
 & \beta=\colP1{k+1} \ar@{^{(}->}[dr]   & \\
\alpha=\colP1k \ar@{^{(}->}[dr] \ar@{^{(}->}[ur] & & \gamma=\colP1{k+2} \\
 & \delta = \{ P_1, \cdots ,P_k, P_{k+2} \} \ar@{^{(}->}[ur] &
}
\end{equation*}
Let $f$ be a discrete Morse function on the simplicial
complex.
Suppose that $f(\alpha) \geq f(\beta) \geq f(\gamma)$, that
is both  \eqref{eq:22} and \eqref{eq:23} hold.
We also have $f(\alpha) < f(\delta) < f(\gamma)$, or
\begin{equation*}
  f(\gamma) > f(\delta) > f(\alpha) \geq f(\beta) \geq f(\gamma)
\end{equation*}
This is a contradiction, so we see that on simplicial complexes
at most one \eqref{eq:22} and \eqref{eq:23} can hold.
\par
However this argument is not restricted to simplicial complexes.
All that is needed is that if $\alpha$ and $\gamma$ are two
cells with $\alpha\subset\gamma$ and 
$2+\dim\alpha=\dim\gamma$ then there are at least two cells in between
$\alpha$ and $\gamma$. 
\begin{theorem}\label{thm:simplicialproperty}
 Between
any two cells  $\alpha$ and $\gamma$ of a Delaunay tesselation with $2+\dim\alpha=\dim\gamma$
and $\alpha\subset\gamma$
there are $\beta_1$ and $\beta_2$, both cells in the Delaunay tesselation, 
both having dimension $1+\dim\alpha$ such that
$\alpha\subset\beta_i\subset\gamma$. 
\end{theorem}
\begin{proof}
Consider the cone $\Cone(\gamma,\alpha)$ of $\gamma$ over $\alpha$:
\begin{equation*}
\Cone(\gamma,\alpha)\overset{\textrm{def}}{=}\{
\beta \,\mid\, \alpha\subseteq \beta\subseteq\gamma\}
\end{equation*}
Intersect that cone with the $n-k+1$ dimensional plane through zero
and orthogonal to the $k-1$-dimensional plane affinely spanned by
$\alpha$. The intersection is a cone in the two dimensional plane.
It has two extremal vectors $\xi_1$ and $\xi_2$. These correspond
to the faces $\beta_1$ and $\beta_2$ we were looking for.
\end{proof}
Hence even if the Delaunay tesselation is not simplicial,
then we still have that not both \eqref{eq:22} and \eqref{eq:23} can
hold. 
\begin{lemma}
For a Delaunay tesselation with a discrete Morse function,
then for any cell $\beta$, either:
\begin{equation}\label{eq:222}
\# \{ \alpha\in\mathcal{T}\,\mid\, 1 + \dim(\alpha)=\dim(\beta)\,\,\,
 \alpha\subset\beta
 \,\,\, h(\alpha)\geq h(\beta) \} = 0
\end{equation}
or 
\begin{equation}\label{eq:233}
\# \{ \alpha\in\mathcal{T}\,\mid\, \dim(\alpha)=1+ \dim(\beta)\,\,\,
 \beta\subset\alpha
 \,\,\, h(\alpha)\leq h(\beta) \} = 0 \ .
\end{equation}
\end{lemma} 
This will play a crucial role in the next section.
\subsection{Discrete vector fields.}
On any coherent triangulation there exists a discrete Morse function:
$h\colon\alpha\mapsto\dim\alpha$.
But this function has all cells critical! 
The construction of discrete Morse functions with less critical points
(even a minimal number), or prescribed values on certain cells is more delicate. 
As an example, in \cite{KKM} one starts with assigning function values to vertices of a simplicial complex.
They construct an algorithm for extension of a function on the vertices
of a simplicial complex to a discrete Morse function on all cells of simplicial complex. 
We  consider an other extension problem. The construction is rather involved and both
\cite{KKM} and we below use the concept of discrete gradient vector field from \cite{MR1612391}.
\par 

Let us have a discrete
Morse function $h$ on the coherent triangulation. The (discrete) gradient of
$h$ consists (by definition) of arrows $\alpha\rightarrow\beta$ between cells of adjacent dimensions, drawn when
$\alpha\subset\beta$ and $h(\alpha)\geq h(\beta)$. 
These cells can
only occur in pairs, which are disjoint. Keeping in mind the values of $h$ we call $\alpha$ {\it up}  and $\beta$ {\it down}. The remaining cells are all critical.
In general such a collection of arrows is called a
\emph{discrete vector field}
A $V$-path is a sequence of arrows such that the cell pointed to 
contains a second 
cell that points to another one, of higher
dimension.
A \emph{closed $V$-path} is a circular $V$-path as in
equation \eqref{eq:10}.
Theorem 3.5 in \cite{MR900443} states that
a discrete vector field is the gradient of a discrete Morse function
if and only if it contains no closed $V$-paths.
\begin{equation}\label{eq:10}
\alpha_1 \rightarrow \beta_1 \supset 
\alpha_2 \rightarrow \beta_2 \supset \cdots \supset
\alpha_k \rightarrow \beta_k \supset \alpha_1
\end{equation}
The dual to the Delaunay tesselation, the power diagram, is also a
polyhedral complex. Given a discrete Morse function $\mathfrak{g}$ on the 
Delaunay tesselation, clearly $-\mathfrak{g}$ is a discrete Morse function
on the power diagram. 
\subsection{Collapses}
In section \ref{sec:morseposetToVf} polyhedral collapses will play a significant role.
To prepare for that section we recall some background information, see \cite{MR1997069}
and \cite{Mi2} for more details.
\par
Consider a polytope $A$ with one of its codimension $1$ faces $F$.
We will delete the interior of $A$ and the interior of the face $F$.
The resulting space $B$ is the union of the remaining codimension $1$ faces of $A$.
We call this process an elementary collapse of the polytope $A$.
Notation $A \searrow B$.
\par
Assume next that we have a space $K$ tessellated by polytopes and that the polytope $A$ is a
$k$-cell of $K$. Assume that $A$ has a free face $F$, i.e.~$F$ is only a face of $A$ and
not a face of any other cell.
We repeat the construction and consider the space $L$ obtained by deleting from $K$
the interior of $A$ and the interior of the face $F$.
The transition from the space $K$ to the space $L$ is called an {\em elementary collapse}.
Notation $K \searrow L$.
A sequence of elementary collapses is called a {\em polyhedral collapse}.
\par
Collapses are examples of deformation retractions in algebraic topology.
Whitehead introduced simple homotopy theory. A basic reference for this is Milnor's
paper \cite{Mi2}. 
Two spaces are called $s$-homotopy equivalent if there is a sequence of elementary
collapses and expansions to go from one space to the other.
It is clear that any $s$-homotopy equivalence determines homotopy equivalence. 
The converse is not true. Whitehead's result states that
there is an obstruction in the so-called Whitehead group
of the fundamental group of the spaces. As soon as this obstruction is zero,
which happens e.g. in simply connected spaces, both equivalences coincide.
\par
Given an elementary collapse we could attach an arrow of a discrete vector field to
the pair of cells.
But also as soon as we have free faces, an arrow of a discrete vector field can be used for
describing an elementary collapse. Note that in this way a $V$-path
\begin{equation*}
\alpha_1 \rightarrow \beta_1 \supset 
\alpha_2 \rightarrow \beta_2 \supset \cdots \supset
\alpha_k \rightarrow \beta_k
\end{equation*}
is naturally related to a s-homotopy equivalence:
\begin{equation*}
L_1 \searrow  L_2 \searrow \cdots \searrow L_k .
\end{equation*}
\section{From Morse poset to a discrete vector field}\label{sec:morseposetToVf}
\subsection{Introduction}
In this section we will construct the discrete Morse function on the Delaunay tesselation
from the power distance function $g$ defined in~\eqref{eq:20}:
\begin{theorem}\label{sec:from-morse-poset}
On $\Del(\colP1N)$  there exists a discrete Morse function $h$ such that the
Morse poset of $h$ equals the Morse poset of $g(x)=\min_{1\leq i \leq N}g_i(x)$.
\end{theorem}
In fact we will not construct $h$ directly. Instead we will construct
a discrete vector field without closed $V$-paths.  
\par
We will use three different types of cells: {\it active}, {\it up} and {\it down}.
The important difference with the gradient vector field above is that up-cells and
down-cells no longer will occur in pairs, but up-cells can be coupled to several down-cells
(of lower dimensions).
\par
In the construction below the power center $c(\alpha)$ will play an important role.
\par
First recall definition \ref{def:active}:
The cell $\alpha$ in the Delaunay tesselation  is \emph{active} or \emph{critical} if
\begin{itemize}
\item[-] $c(\alpha)$ is contained in the relative interior
$\Interior( \CH( \alpha ) )$ of $\CH(\alpha)$, and
\item[-] $g( c(\alpha) ) = g_i(c(\alpha))$ for all $P_i \in \alpha$ and
 $g( c(\alpha) ) < g_j(c(\alpha))$ for all $P_j \not \in \alpha$.
\end{itemize}
\par
A face $\alpha$ can be non-active for one of two reasons
\begin{itemize}
\item[] \underline{Down}: The separator $\Sep(\alpha)$ can lie outside
 $\Interior(\CH(\alpha))$, or if that is not the case:
\item[] \underline{Up}: Some point $P_i\in\colP1N\setminus\alpha$ can lie closer to
$c(\alpha)$ than the points in $\alpha$, i.e.\
\begin{equation}\label{eq:5}
  g_i(c(\alpha))\leq g_j(c(\alpha))\quad \forall j\in \alpha
\end{equation}
\end{itemize}
NB. We have identified $\alpha$ with its index set.
\par
To illustrate these concepts we refer the reader to figure \ref{fig:dried}.
In that figure the cell $\P12$ is not active for the ``Up'' reason.
The cell $\PPP1234$ is not active. It cannot be not active for the
``Up'' reason, because there is no cell in which it is contained. Indeed
$\PPP1234$ is not active for the ``Down'' reason. The separator 
$\Sep(P_1, P_2, P_3 , P_4)$ is just $c(P_1, P_2, P_3 , P_4)$ and it clearly
lies outside the tetrahedron.
\par
The reader is encouraged to repeat this
exercise for figure \ref{fig:13}, which we have also drawn in figure \ref{fig:14}.
\par
Next we show two properties:
\begin{itemize}
\item[-]
Suppose that $\alpha$ is not active for the ``Down'' reason.
Then in the case of Voronoi diagrams $\alpha$ cannot be a vertex or an edge.
In general with power diagrams $\alpha$ cannot be a vertex.
\item[-]
Suppose that $\alpha$ is not active for the ``Up'' reason and that 
$\dim\Del(\alpha)=n$. Then we have 
\begin{equation*}
  \{c(\alpha)\}= \cap_{j\in\alpha}\Pow(\{P_j\})\text{:}
\end{equation*}
we cannot have inequality \eqref{eq:5}. We conclude that no $\alpha$, not active for the
``Up'' reason has $\dim\Del(\alpha)=n$.
\end{itemize}
\subsection{Investigating the up-down structure}
We start by showing that for every face, not active for the 
``Down'' reason, there is a unique lower dimensional cell
not active for the ``Up'' reason. Then we construct for each
face not active for the ``Up'' reason a part of the discrete
vector field. 
\begin{lemma}\label{thm:downreason}
Let $\beta\in\Del(\colP1N)$ be a face not active
for the ``Down'' reason. Then the closest point to
$\CH(\beta)$ from $c(\beta)$ is $c(\alpha)$ for
some proper face $\alpha$ of $\beta$. At $c(\alpha)$ we have
\begin{equation}\label{eq:6}
  g_j(c(\alpha)) \leq g_\alpha(c(\alpha)) \quad \forall j\in\beta\setminus\alpha
\end{equation}
Hence $\alpha$ is not active for the ``Up'' reason.
Conversely if \eqref{eq:6} holds for $\beta\supset\alpha$ then $\beta$ is not
active for the ``Down'' reason.
\end{lemma}
\begin{proof}
We start with ``$\Rightarrow$''.
\par
A special case we first have to handle is when 
$c(\beta)$ lies on the relative boundary of
$\CH(\beta)$. This means that $c(\beta)\in\CH(\alpha)$ for
some proper face $\alpha$ of $\beta$. Because $\Sep(\beta)\subset
\Sep(\alpha)$ it follows that $c(\beta)=c(\alpha)$. Equation
\eqref{eq:6} follows automatically.
\par
Now we can safely assume that the distance from $c(\beta)$ to
$\CH(\beta)$ is $>0$.
Denote $y$ the closest point on $\CH(\beta)$ from $c(\beta)$.
If $\alpha$ is a vertex we are done, so suppose that $\alpha$
is not a vertex. The line
from $y$ to $c(\beta)$ is orthogonal to $\Aff(\alpha)$ and
hence parallel to $\Sep(\alpha)$. Clearly $c(\beta)$ lies
in $\Sep(\alpha)$. So $y$ also lies in $\Sep(\alpha)$. 
But $y$ also lies in $\Aff(\alpha)$. So, by the definition
of $c(\alpha)$ in equation \eqref{eq:24} we have that
$y=c(\alpha)$.
\par
Because $\CH(\beta)$ is a convex set the hyperplane with
normal $c(\beta)-c(\alpha)$ passing through $c(\alpha)$ separates
$\Interior(CH(\beta))$ from $c(\beta)$. 
Thus for a vertex $P_k$ of $\beta$ the angle $P_k$ to $c(\alpha)$ to $c(\beta)$ must be obtuse.
Consequently:
\begin{equation}
  \label{eq:7}
  \lVert c(\beta) - P_k \rVert^2 \geq \lVert c(\beta) - c(\alpha) \rVert^2
+ \lVert c(\alpha) - P_k \rVert^2
\end{equation}
if $P_k$ is a vertex in $\beta\supset\alpha$.
By the Pythagoras theorem, equality holds if $k\in\alpha$:
\begin{equation}
  \label{eq:9}
  \lVert c(\beta) - P_k \rVert^2 = \lVert c(\beta) - c(\alpha) \rVert^2
+ \lVert c(\alpha) - P_k \rVert^2
\end{equation}
Then equation \eqref{eq:7} becomes, taking the factor $\frac12$ we 
used in the definition of $g$ in equation \eqref{eq:20} into account:
\begin{equation}
  \label{eq:8}
  g_k(c(\beta)) \geq \frac12 \lVert c(\beta) - c(\alpha) \rVert^2 + g_k(c(\alpha))
\quad k \in \beta\setminus\alpha
\end{equation}
Again, equality holds when $k\in\alpha$.
Putting this together we get \eqref{eq:6}, and we see that $\alpha$ is not
active for the ``Up'' reason.
\par
Next do the ``$\Leftarrow$'' part. We assume \eqref{eq:6} and $\alpha\subset\beta$.
From \eqref{eq:6} we get for $j\in\beta\setminus\alpha$:
\begin{equation}\label{eq:1}
  g_j(c(\alpha)) +\frac12 \lVert c(\beta) - c(\alpha) \rVert^2 \leq
g_\alpha(c(\alpha))+\frac12 \lVert c(\beta) - c(\alpha) \rVert^2 = g_\alpha(c(\beta))=
g_j(c(\beta))
\end{equation}
The special case to handle first is where equality holds in 
\eqref{eq:1}. Then $c(\beta)=c(\alpha)$ and so $c(\beta)$ is
not an element of the relative interior $\Interior(\CH(\beta))$.
So $\beta$ is not active for the ``Down'' reason.
\par
Now we can safely assume strict inequality in \eqref{eq:1}.
The simplest case is a $\beta$ for which $\beta = \alpha \cup \{ j \}$.
The three points $P_j$, $c(\beta)$ and $c(\alpha)$ all lie in $\Aff(\beta)$.
From \eqref{eq:1} we get \eqref{eq:7} with strict inequality. 
The line segment from $c(\alpha)$ to $c(\beta)$ is thus an outward
normal to $\CH(\beta)$ in $\Aff(\beta)$. Consequently $\beta$ is not active for the
``Down'' reason. The general case is similar.
\end{proof}
\subsection{Defining the up-sets}
Let again $\alpha$ be not active for the ``Up'' reason.
Denote $K$ the set of indices for which \eqref{eq:5} holds.
Denote $\Up(\alpha)$ the set of simplices:
\begin{equation}\label{eq:4}
  \Up(\alpha)=\{ \beta\in\Del(\colP1N) \,\mid\,
\alpha\subset\beta\subset\alpha\cup K \}
\end{equation}
Lemma \ref{thm:downreason} characterizes the elements of $\Up(\alpha)$
as those $\beta\supset\alpha$ for which $c(\alpha)$ is the closest point on
$\CH(\beta)$ from $c(\beta)$. The complement of the Morse poset is 
divided into different subsets $\Up(\alpha)$, one for each face not active
for the ``Up'' reason. 
We need to prove that each $\Up(\alpha)$ can be filled up
with a discrete vector field. 
\par
\begin{figure}[htbp]
  \centering
  \includegraphics[height=2.5cm]{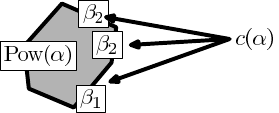}
  \caption{View from $c(\alpha)$ to $\Pow(\alpha)$ inside the plane $\Aff(\alpha)$.%
In this case $\Up(\alpha)$ contains $\beta_1$, $\beta_2$ and $\beta_3$. Lemma \ref{thm:up-polyhedron} shows that this is equivalent to their power cells being ``visible'' from $c(\alpha)$. }
  \label{fig:6}
\end{figure}
\begin{figure}[htbp]
  \centering
  \includegraphics[height=6cm]{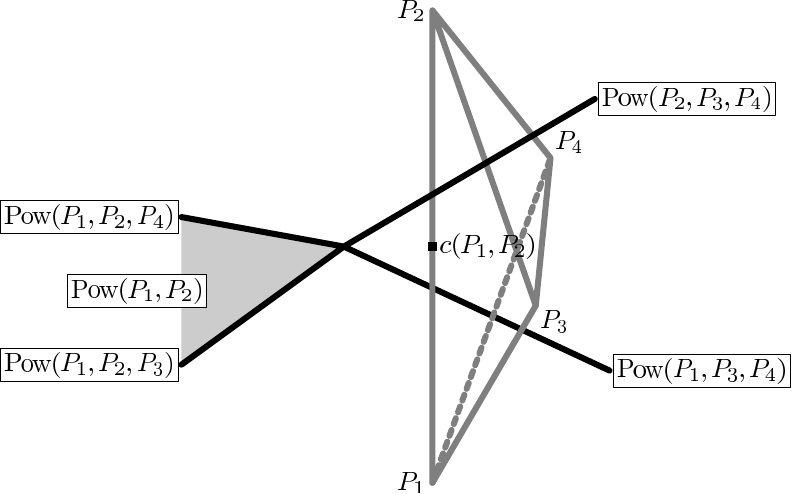}
  \caption{Second, less schematic example, of $\Up(\alpha)$. In this example $\Up(P_1,P_2)$ consists of $\P12$, $\PP123$, $\PP124$ and $\PPP1234$. From $c(P_1,P_2)$ we can see exactly the power cells of those simplices that lie in $\Up(P_1,P_2)$.}
  \label{fig:dried}
\end{figure}

We will establish another criterion (in terms of power sets) for which a $\beta\in\Del(\colP1N)$ with
$\alpha\subset\beta$ is an element of $\Up(\alpha)$.
Take the point $c(\alpha)$ outside the polytope $\Pow(\alpha)$ in $\Sep(\alpha)$
and consider the faces of $\Pow(\alpha)$ for which
the outward pointing normal makes an angle smaller
than 90 degrees with $x-c(\alpha)$, where $x$ be a point in the interior of $\Pow(\beta)$. Then we get
a number of faces $\beta_1, \cdots , \beta_l$ in
the Delaunay tesselation that we see as faces $\Pow(\beta_1), \cdots , 
\Pow(\beta_l)$ on the relative boundary $\partial(\Pow(\alpha))$ of $\Pow(\alpha)$.
This is illustrated in figure \ref{fig:6} and \ref{fig:dried}. 
A precise formulation is contained in the following lemma.
\begin{lemma}\label{thm:up-polyhedron}

Let $\beta\supset\alpha$ be an element of $\Del(\colP1N)$, $\alpha$ not active for the ``Up'' reason.
Let $x$ be a point in the interior of $\Pow(\beta)$.
If, either $x=c(\alpha)$, or 
\begin{enumerate}
\item the line through $x$ and $c(\alpha)$ intersects $\Pow(\alpha)$ in a
segment, and 
\item the segment $[x,c(\alpha)]$ contains no point of the interior of
$\Pow(\alpha)$
\end{enumerate}
then $\beta\in\Up(\alpha)$.
Conversely, if $\beta\in\Up(\alpha)$ and $x$ lies in the interior of 
$\Pow(\beta)$ then 1 and 2 hold, or $x=c(\alpha)$.
\end{lemma}
\begin{proof}
We first need to treat the case $x=c(\alpha)$.
In that case $c(\alpha)$ lies on the relative boundary of $\Pow(\alpha)$.
So $x\in \Pow(\beta)$, $\beta=\alpha\cup I$.
And $f_j(c(\alpha)) = f_\alpha(c(\alpha))$, for all $j\in I$. 
(note that we use here the affine functions from section \ref{sec-power}.)
Obviously $\Up(\alpha)$ consists exactly of all those $\gamma\in\Del(\colP1N)$
with
\begin{equation*}
\alpha\subseteq \gamma \subseteq \beta=\alpha\cup I
\end{equation*}
\par
Now we can assume $x\not=c(\alpha)$. Again $x$ lies on the relative
boundary of $\Pow(\alpha)$ and $\beta=\alpha\cup I$ for some maximal
set $I$.
We put $x_t=tx+(1-t)c(\alpha)$. Hence $x_1=x$ and $x_0=c(\alpha)$.
It also follows for $j\in I$
\begin{equation*}
  f_j(x_t)= t f_j(x)+(1-t)f_j(c(\alpha))=
t f_\alpha(x) +(1-t)f_j(c(\alpha))
\end{equation*}
and thus
\begin{equation*}
  f_j(x_t)-f_\alpha(x_t)=(1-t)(f_j(c(\alpha))-f_\alpha(c(\alpha)))
\end{equation*}
Now we have
\begin{equation*}
f_j(c(\alpha))\geq f_\alpha(c(\alpha))\,\Leftrightarrow\,
\exists t>1\,\colon\, f_j(x_t)\leq f_\alpha(x_t)
\end{equation*}
So that
\begin{equation*}
\exists t > 1\,\, x_t\in\Pow(\alpha)\,
\Rightarrow\,f_j(c(\alpha))\geq f_\alpha(c(\alpha))
\end{equation*}
And this is exactly what we needed to prove.
\par
For the converse suppose that $\beta\supset\alpha$, and $\alpha$ not active
for the ``Up'' reason, but that 1 and 2 do not hold.
\par
The only point of $\Pow(\alpha)$
on the line through $x$ and $c(\alpha)$ is $x$.
It follows that $\beta\setminus\alpha$ contains
more than one index, i.e.\ $\beta=\alpha\cup\{j_1 , \cdots , j_r\}$ with
$r\geq 2$. In a sufficiently small neighborhood of $x$ the half spaces
$H_i=\{f_{j_i}\leq f_\alpha\}$ define $\Pow(\alpha)$ as a polyhedron
in $\Sep(\alpha)$. Because $\{x_t\,\mid\,t\in\RR\}$ touches $\Pow(\alpha)$
there exist an index $k$, say $k=j_1$, and a $t<1$ such that
$f_j(x_t)-f_\alpha(x_t)$ and hence $f_k(c(\alpha))<f_\alpha(c(\alpha))$.
And so $\beta\supset\alpha$ is not an element of $\Up(\alpha)$.
\end{proof}

\subsection{Construction of the discrete vector field}
Let $\alpha$ be not active for the ``Up'' reason.  Hence $c(\alpha)$ lie outside $\Pow(\alpha)$.
We will consider $c(\alpha)$ as center of a system of rays. As seen above this system of rays defines a cone, which meets the boundary $\partial \Pow(\alpha)$ exactly in $ \text{Front} (\alpha)$, i.e. those $\Pow(\beta)$, where $\beta$ in $\Up(\alpha)$ is different from $\alpha$.
\par
This $ \text{Front} (\alpha)$ is part of  $\partial \Pow(\alpha)$ and is contractible, since $\Pow(\alpha)$ is convex.
\par
We can use the extension of rays to define a deformation retract to the back side:
$$ 
\Pow(\alpha)   \longrightarrow \text{Back}(\alpha) := \overline{ \Pow(\alpha) - \text{Front} (\alpha) }
$$
NB.\ In case $\Pow(\alpha)$ is not bounded there is  a small complication, which can be solved by adding some extra points, that do not effect the situation in $\Up(\alpha)$.
\par
Since all spaces we consider are contractible, this deformation retraction can be realized
by elementary collapses. This a different way to express, that there exists a discrete vector
field on $\Pow ((\Up (\alpha))$. Due to the contractibility this vector field has no closed
V-paths.
\par
By duality there exists a discrete vector field without closed V-paths on $\Up(\alpha)$.
\par
We conclude that the complement of the Morse poset
can be divided into separate parts $\Up(\alpha)$, for all
up-cells $\alpha$ and on each $\Up(\alpha)$
we have defined a discrete vector field without closed V-paths.
This defines a discrete vector field on $\Del(P_1,\cdots,P_N)$ by ``no-arrow'' extension on the Morse poset.
To show that we indeed get a discrete Morse function  we
need to prove one more thing: there should be no global closed
$V$-path.
\subsection{No global V-loops}
\begin{figure}[htbp]
  \centering
  \includegraphics[height=2cm]{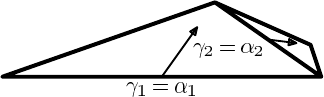}
  \caption{A $V$-path, with a jump from $\Up(\gamma_1)$ to $\Up(\gamma_2)$.}
  \label{fig:jump}
\end{figure}

Suppose that we have a closed $V$-path. Inside one $\Up(\alpha)$ there is
no such closed $V$-path, so at some cell the $V$-path should jump
from one $\Up(\gamma_1)$ to another $\Up(\gamma_2)$, as in figure
\ref{fig:jump}.
\begin{lemma}\label{sec:from-morse-poset-4}
Suppose that we have a $V$-path that jumps from $\Up(\gamma_1)$ to $\Up(\gamma_2)$
as in equation \eqref{eq:3}.
\begin{equation}\label{eq:3}
  \begin{matrix}
\alpha_i & \rightarrow & \beta_i & & & \subset & \Cone(\gamma_1\cup I_1,\gamma_1) \\
 & & \cup & & & & \\
 & &\alpha_{i+1} & \rightarrow & \beta_{i+1} & \subset & \Cone(\gamma_2\cup I_2,\gamma_2)
  \end{matrix}
\end{equation}
Then $g_{\gamma_1}(c(\gamma_1))>g_{\gamma_2}(c(\gamma_2))$.
\end{lemma}
\begin{proof}
Because $\gamma_1$ is a proper subset of
$\gamma_1\cup\gamma_2$ we see that:
$c(\gamma_1\cup\gamma_2)\not\in\Interior(\CH(\gamma_1\cup\gamma_2))$.
The closest point to $\CH(\gamma_1\cup\gamma_2)$ from $c(\gamma_1\cup\gamma_2)$ is
$c(\gamma_1)$, as we have shown in lemma \ref{thm:downreason}. In particular
\begin{equation}\label{eq:11}
  \lVert c(\gamma_1 \cup \gamma_2) - c(\gamma_1) \rVert^2 <
  \lVert c(\gamma_1 \cup \gamma_2) - c(\gamma_2) \rVert^2 
\end{equation}
Here, the inequality is strict. The closest point is unique, because
$\CH(\gamma_1\cup\gamma_2)$ is convex.
\par
For the reasons explained in the proof of lemma \ref{thm:downreason}
we have
\begin{equation}\label{eq:12}
  g_{\gamma_1\cup\gamma_2}(c(\gamma_1\cup\gamma_2))=
  g_{\gamma_1}(c(\gamma_1\cup\gamma_2))=
\frac12 \lVert c(\gamma_1\cup\gamma_2) - c(\gamma_1) \rVert^2 + g_{\gamma_1}(c(\gamma_1))
\end{equation}
The same identity holds with $\gamma_1$ and $\gamma_2$ exchanged, thus:
\begin{equation}
  \label{eq:13}
\frac12 \lVert c(\gamma_1\cup\gamma_2) - c(\gamma_1) \rVert^2 +
 g_{\gamma_1}(c(\gamma_1)) =
\frac12 \lVert c(\gamma_1\cup\gamma_2) - c(\gamma_2) \rVert^2 +
 g_{\gamma_2}(c(\gamma_2))
\end{equation}
Putting \eqref{eq:11} and \eqref{eq:13} together we get the
desired result.
\end{proof}
Lemma \ref{sec:from-morse-poset-4} says that to each $\Up(\alpha)$ appearing
in the $V$-path we can associate a number, and
passing from one $\Up(\alpha)$ to another that number strictly decreases. 
Hence we can not have a closed $V$-path. The proof of theorem
\ref{sec:from-morse-poset} is complete.

\subsection{Constructing the extension}
We constructed above a discrete Morse function with the same Morse poset as $g$.
What is lost in the proof is a direct relation between $g$ and the 
discrete Morse function that results from the discrete vector field
we constructed.
A good candidate for such a Morse function could be:
\begin{equation*}\label{eq:14}
h: \mathfrak{g} \to  \RR  \ \  ; \ \ \    h(\beta) = \sup_{x\in\CH(\beta}g(x)
\end{equation*}
as a function on the Delaunay tesselation.
Inspection of this function shows, that it is indeed critical on the Morse poset, 
but has constant values on all cells of a given $\Up(\alpha)$.
So it violates the rules of Morse functions as soon as $\Up(\alpha)$ contains more than two cells.
We will use the discrete vector field on $\Up(\alpha)$ (defined above)  to perturb these values.
Remark that the vector field induces on $\Up(\alpha)$ a `potential function'; call it $h_{\alpha}$.
Set $h_{\alpha}(\beta) = 0$ as soon as $\beta \not \in \Up(\alpha)$.
One has the freedom to scale $h_{\alpha}$ such that it takes values in small interval around 0.
The function
\begin{equation*}\label{eq:18}
\hat{h}: \mathfrak{g} \to  \RR  \ \  \ ;  \ \    \hat{h}(\beta) = 
\sup_{x\in\CH(\beta)}g(x) + \sum_{\alpha \ \text{up}} h_{\alpha}(\beta)
\end{equation*}
is a good extension as soon as the scaling is small enough. So we have :
\begin{theorem}\label{sec:morse-ext}
The function $\hat{h}$ is a discrete Morse function, with the same Morse poset as $g$ and which assigns to each element of the Morse poset the corresponding critical value of $g$.
\end{theorem}
\begin{remark}\label{sec:from-morse-poset-2}
The above discussion shows that the vector field is not unique. The non-uniqueness only lies
in the choice of the arrows possible in $\Up(\alpha)$.
It seems more natural to consider a single polyhedral collapse of $\Up(\alpha)$.
But Forman's discrete Morse theory is ill equipped for dealing with bigger polyhedral
collapses. 
One could also try to develop a generalization of discrete Morse theory
in order to include these polyhedral collapses.
\end{remark}
\section{Examples and comments.}\label{sec-examples}
\subsection{Vertices that are not active}
\begin{figure}[htbp]
  \centering
  \includegraphics[height=3cm]{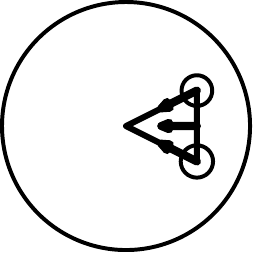}
  \caption{Morse poset and discrete vector field of figure \ref{fig:13}}
  \label{fig:14}
\end{figure}

Let us return to the counterintuitive example of figure 
\ref{fig:13}. There the Morse poset is as in figure
\ref{fig:14}. In the picture we see three cones,
consisting each of two faces in the simplest triangulation
one can think of. 
\par
\subsection{3 dimensional case}
In \cite{SieVanM} we studied the Morse theory of a tetrahedron in $\RR^3$.
As we found, that for triangles in $\RR^2$ generically there are two different Morse posets,
we showed that there are generically nine different Morse posets for the tetrahedron.
To be more precise: theorem \ref{thm:tetra}
tells us exactly the possible positions of critical points.
The number $s_i$ of critical point of index $i$ are bounded by $4,6,4,1$ for index $0,1,2,3$
Moreover the Euler characteristic of $\RR^3$ is $+1$, so by Morse Formula: $s_0-s_1+s_2-s_3=1$
This gives a priori the following $9$ possibilities:
\begin{equation*}
\begin{array}{|c|c|c|c|}
\hline
s_{0} & s_{1} & s_{2} & s_{3} \\
\hline
4 & 6  & 4  & 1 \\
4 & 5  & 3  & 1 \\
4 & 4  & 2  & 1 \\
4 & 3  & 1  & 1 \\
4 & 2  & 0  & 1 \\
4 & 6  & 3  & 0 \\
4 & 5  & 2  & 0 \\
4 & 4  & 1  & 0 \\
4 & 3  & 0  & 0 \\
\hline
\end{array}
\end{equation*}
But not all these will occur.
Since we start with $4$ vertices and the result should be a connected space,
we need at least $3$ saddle points of index $1$.
Hence, we cannot have $(4,2,0,1)$.
\par
\begin{figure}[hbtp]
\centering
\includegraphics[height=9cm]{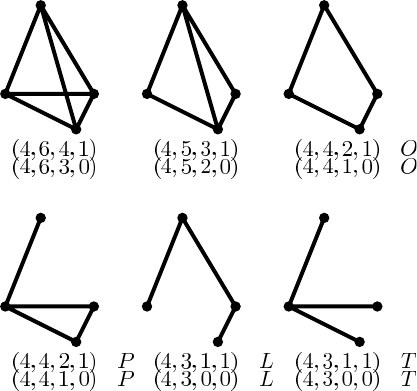}
\caption{List of a priori possible Gabriel graphs. The vertices are the minima of the distance function. 
Each vertex of the tetrahedron $\mathbb{T}$ is a minimum of $d$. There are no other minima. 
The edges of the graph are the index $1$ critical points of the distance function. 
Not each midpoint of an edge of $\mathbb{T}$ is an index $1$ critical point of the distance function $d$. 
The graph $(4410)\, O$ can be laid out so as to form the letter ``O''. The names of the other graphs are chosen similarly.}
\label{fig:1}
\end{figure}

We list in figure \ref{fig:1} the (a priori) possible 1-skeletons of the Morse poset
(also know as Gabriel graphs) for the above cases, they
are the connected graphs with $4$ vertices.
Just as the case $(4,2,0,1)$ can not occur, we showed
in \cite{SieVanM}
 that the cases $(4,4,2,1)\, P$,
$(4,3,1,1)\, T$ and $(4,3,1,1)\, L$ do not occur. 
\begin{theorem}\label{thm:tetra}
Up to affine bijections of a tetrahedron in $\RR^3$, 
which send the Morse posets to each other, 
there are nine generic tetrahedra. They are uniquely described by the 
nine Gabriel graphs $(4,6,4,1)$, $(4,6,3,0)$, $(4,5,3,1)$,
$(4,5,2,0)$, $(4,4,2,1)\, O$, $(4,4,1,0)\, O$, $(4,4,1,0)\, P$,
$(4,3,0,0)\, L$ and $(4,3,0,0)\, T$, drawn in figure
\ref{fig:1}. 
\end{theorem}
The example $4300I$ is interesting, since it shows a polyhedral collapse.
In figure \ref{fig:dried} we have reproduced an example of such 
a tetrahedron.
The face $\P12$ is not active for the ``Up'' reason.
In that case both $\PP123$ and $\PP124$ and $\PPP1234$ are not
active for the ``Down'' reason. There is one cone
consisting of the polyhedral collapse:
\begin{equation*}
  \P14 \rightarrow \PP124 \quad \PP134 \rightarrow \PPP1234 
\end{equation*}
The decomposition of that cone into arrows is not unique.
We might just as well write
\begin{equation*}
  \P14 \rightarrow \PP134 \quad \PP124 \rightarrow \PPP1234 
\end{equation*}
\subsection{2 dimensional case}
\begin{figure}[htbp]
  \centering
  \includegraphics[width=\textwidth]{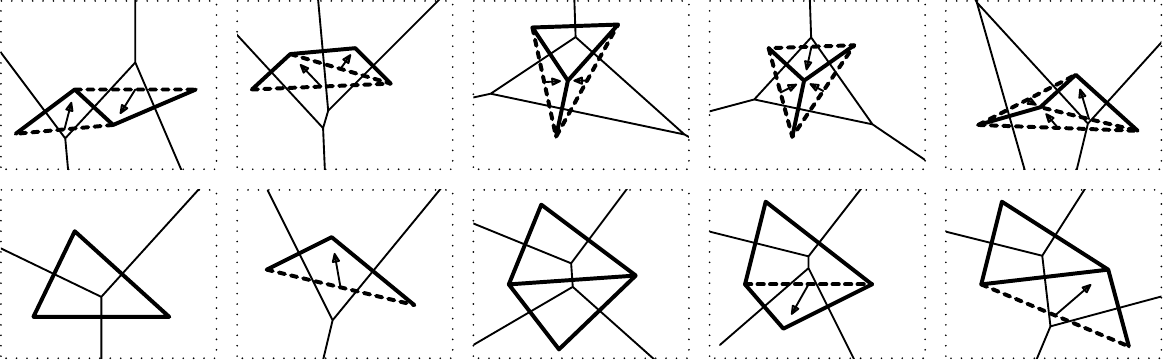}
  \caption{$3$ and $4$ points in $\RR^2$: all Morse posets from Voronoi diagrams with the accompanying discrete vector fields}
  \label{fig:8}
\end{figure}

In the same spirit we have the classification of all Morse posets
of Voronoi diagrams for $4$ points in the plane in \cite{VoronoiDist}.
Unfortunately there were two cases missing, which are now included in figure
\ref{fig:8}. 
\subsection{Enumerating all power diagrams and Voronoi diagrams}
The Morse poset is the discrete structure that encodes the flow of the
distance function $g(x)=\min_{1\leq i \leq N} g_i(x)$.
It is unclear whether the problem of enumerating
all power diagrams is equivalent to the problem of enumerating all coherent
triangulations. That last problem is nicely
solved using the secondary polytope, see \cite{MR1264417}.
For many applications
though a discrete structure describing a Voronoi diagram should
contain more geometric information. The Morse poset and the discrete vector field
do exactly that. In figure \ref{fig:8} we see that there are 2 combinatorially
different triangulations of four points in the plane. Taking into
account the Morse poset yields 8 different cases.
\bibliographystyle{plain}

\end{document}